\documentclass{amsart}
\usepackage{amsmath} 
\usepackage{amssymb}
\usepackage{mathrsfs}

\newtheorem{theorem}{Theorem}[section] 
\newtheorem{claim}[theorem]{Claim}

\newtheorem{conclusion}[theorem]{Conclusion}

\theoremstyle{definition}
\newtheorem{definition}[theorem]{Definition}
\newtheorem{example}[theorem]{Example}
\newtheorem{thesis}[theorem]{Thesis}
\newtheorem{conjecture}[theorem]{Conjecture}
\newtheorem{discussion}[theorem]{Discussion}

\theoremstyle{remark}
\newtheorem{remark}[theorem]{Remark}
\newtheorem{question}[theorem]{Question}
\newtheorem{notation}[theorem]{Notation}

\newcommand{\rest}{{\restriction}}

\newcommand{\tr}{{\rm tr}}
\newcommand{\wilog}{{\rm without loss of generality}}

\newcommand{\then}{{\underline{then}}}
\newcommand{\when}{{\underline{when}}}
\newcommand{\Then}{{\underline{Then}}}

\newcommand{\If}{{\underline{if}}}
\newcommand{\Iff}{{\underline{iff}}}
\newcommand{\mn}{{\medskip\noindent}}
\newcommand{\sn}{{\smallskip\noindent}}
\newcommand{\bn}{{\bigskip\noindent}}

\newcommand{\cA}{{\mathscr A}}

\newcommand{\bbF}{{\mathbb F}}

\newcommand{\gC}{{\mathfrak C}}
\newcommand{\gS}{{\mathfrak S}}

\newcommand{\cR}{{\mathscr R}}

\newcommand{\cF}{{\mathscr F}}
\newcommand{\cG}{{\mathscr G}}
\newcommand{\cI}{{\mathscr I}}

\newcommand{\bbL}{{\mathbb L}}
\newcommand{\cM}{{\mathscr M}}

\newcommand{\bbP}{{\mathbb P}}
\newcommand{\bbZ}{{\mathbb Z}}

\newcommand{\bbQ}{{\mathbb Q}}

\newcommand{\cT}{{\mathscr T}}

\newcount\skewfactor
\def\mathunderaccent#1#2 {\let\theaccent#1\skewfactor#2
\mathpalette\putaccentunder}
\def\putaccentunder#1#2{\oalign{$#1#2$\crcr\hidewidth
\vbox to.2ex{\hbox{$#1\skew\skewfactor\theaccent{}$}\vss}\hidewidth}}
\def\name{\mathunderaccent\tilde-3 }

\newenvironment{PROOF}[2][\proofname.]
   {\begin{proof}[#1]\renewcommand{\qedsymbol}{\textsquare$_{\rm #2}$}}
   {\end{proof}}

\begin{document}

\title {Theories with EF-Equivalent Non-Isomorphic Models}
\author {Saharon Shelah}
\address{Einstein Institute of Mathematics\\
Edmond J. Safra Campus, Givat Ram\\
The Hebrew University of Jerusalem\\
Jerusalem, 91904, Israel\\
 and \\
 Department of Mathematics\\
 Hill Center - Busch Campus \\ 
 Rutgers, The State University of New Jersey \\
 110 Frelinghuysen Road \\
 Piscataway, NJ 08854-8019 USA}
\email{shelah@math.huji.ac.il}
\urladdr{http://shelah.logic.at}
\thanks{The author would like to thank the Israel Science Foundation for
partial support of this research (Grant No.242/03).
I would like to thank Alice Leonhardt for the beautiful typing. Paper Sh:897}

\subjclass{Primary; Secondary:}



\date{December 17, 2012}

\begin{abstract}
Our ``long term and large scale" aim 
is to characterize the first order theories
$T$ (at least the countable ones) such that: for every ordinal
$\alpha$ there are
$\lambda,M_1,M_2$ such that $M_1,M_2$ are non-isomorphic models of $T$
of cardinality $\lambda$ which are EF$^+_{\alpha,\lambda}$-equivalent.
We expect that as in the main gap  we get a strong
dichotomy, so in the non-structure side we have stronger, better examples,
and in the structure side we have a parallel of \cite[Ch.XIII]{Sh:c}.  We
presently prove the \underline{consistency} of the non-structure side for $T$
which is $\aleph_0$-independent (= not strongly dependent), 
even for PC$(T_1,T)$.
\end{abstract}

\maketitle
\numberwithin{equation}{section}
\setcounter{section}{-1}

\centerline {Anotated Content}
\bigskip

\noindent
\S0 \quad Introduction
\bigskip

\noindent
\S1 \quad Games, equivalences and the question
\mn
\begin{enumerate}
\item[${{}}$]   [We discuss what are the hopeful conjectures
concerning versions of EF-equivalent non-isomorphic models for a
given complete first order $T$, i.e. how it fits classification.  In
particular we
define when $M_1,M_2$ are EF$_{\gamma,\lambda}$-equivalent and when
they are EF$^+_{\gamma,\theta,\lambda}$-equivalent and discuss those notions.]
\end{enumerate}
\bigskip

\noindent
\S2 \quad The properties of $T$ and relevant indiscernibility
\mn
\begin{enumerate}
\item[${{}}$]   [We recall the definitions of ``$T$ strongly dependent",
``$T$ is strongly$_4$ dependent", and prove
the existence of models of such $T$ suitable for proving non-structure
theory.]
\end{enumerate}
\bigskip

\noindent
\S3 \quad Forcing EF$^+$-equivalent non-isomorphic models
\mn
\begin{enumerate}
\item[${{}}$]   [We force such an example.]
\end{enumerate}
\bigskip

\noindent
\S4 \quad Theories with order
\mn
\begin{enumerate}
\item[${{}}$]   [We prove in ZFC, that for $\lambda$ regular there
are quite equivalent non-isomorphic models of cardinality $\lambda^+$.]
\end{enumerate}
\newpage

\section {Introduction} 

\subsection {Motivation} \

We first give some an introduction for non-model theorists.  A major
theme in the author's work in model theory is to find ``main gap
theorems".  This means, considering the family of elementary classes
(e.g. the classes of the form Mod$_T =$ the class of models of a
(complete) first order theory $T$, each such class is either very
``simple" or is very complicated; expecting that we have much
knowledge to gain on the ``very simple" ones and even on approximations
to them.

Of course, this depends on the criterion for ``simple".  Essentially
the main theorem of \cite{Sh:c} does this for countable $T$, with
``complicated" interpreted as ``$\dot I(\lambda,T)$", the number of
models in Mod$_T$ of cardinality $\lambda$, is maximal,
i.e. $2^\lambda$, for every $\lambda$.  See more e.g. in
\cite{Sh:E53}.  Here we are interested with interpreting
``complicated"  as ``for arbitrarily large cardinals, there are models
$M_1,M_2 \in \text{ Mod}_T$ of cardinality $\lambda$ which are ``very
similar" but not isomorphic", where ``very similar" is interpreted as
a relative of the game of the following form.  The isomorphism player
constructs during the play, partial isomorphism of cardinality $< \lambda$,
 in each move the anti-isomorphism player demands some
elements to be in the domain or the range, the isomorphism player has to
extend the partial isomorphism accordingly; in the play there are
$\alpha$ moves, $\alpha < \lambda$; and the isomorphism player wins the
play if he has a legal move in each stage (see Definition
\ref{1a.7}, \ref{1a.9}).

In the present paper we try to deal with suggesting the ``right" variant
of the game, (see Definition \ref{1a.8}), and give quite weak sufficient
conditions for Mod$_T$ being complicated.

\centerline {$* \qquad * \qquad *$}

Our aim is to prove (on $\text{PC}(T_1,T)$, see Definition \ref{0z.7}(3))
\mn
\begin{enumerate}
\item[$\boxtimes$]   if $T\subseteq T_1$ are complete first
order theories such that $T$ is not strongly stable, $\alpha$ is
an ordinal and $\lambda > |T|$ (or at least for many such $\lambda$'s)  
\then
\begin{enumerate}
\item[$(*)$]    there are $M_1,M_2 \in \text{ PC}(T_1,T)$ of 
cardinality $\lambda$ which are EF$^+_{\alpha,\lambda}$-equivalent for
every $\alpha < \lambda$ but not isomorphic (for the definition of
EF$^+_{\alpha,\lambda}$, see Definition \ref{1a.9} below, it is a 
somewhat stronger relative of EF$_{\alpha,\lambda}$-equivalent).
\end{enumerate}
\end{enumerate}
\mn

\subsection {Related Works} \

Concerning constructing non-isomorphic
EF$^+_{\alpha,\lambda}$-equivalent 
models $M_1,M_2$ (with no relation to $T$) we
have intended to continue \cite{Sh:836}, or 
see more Havlin-Shelah \cite{HvSh:866} and see history in Vaananen in
\cite{Va95}.  Those works leave the case $\lambda = \aleph_1$ open; 
a recent construction \cite{Sh:907} resolve this but whereas it
applies to every regular uncountable $\lambda$, it seems less amenable
to generalizations. 

By \cite{Sh:c} we essentially know for $T$ a countable complete first
order \when \, there are $\bbL_{\infty,\lambda}(\tau_T)$-equivalent
non-isomorphic models of $T$ of cardinality $\lambda$ for some
$\lambda$, see \S4; this is exactly when $T$ is superstable with NDOP,
NOTOP; (see \cite{Sh:220}).

On restricting ourselves to models of $T$ for
``EF$_{\alpha,\lambda}$-equivalent non-isomorphic", Hyttinen-Tuuri
\cite{HyTu91} started, then Hyttinen-Shelah \cite{HySh:474}, \cite{HySh:529},
\cite{HySh:602}.  
The notion ``EF$^+_{\alpha,\lambda}$-equivalent" is
introduced here, in Definition \ref{1a.9}.  

By \cite{HySh:474}, if $T$ is stable unsuperstable, complete first
order theory, $\lambda = \mu^+,\mu = \text{ cf}(\mu) \ge |T|$,
\then \, there are EF$_{\mu \times \omega,\lambda}$-equivalent
non-isomorphic models of $T$ (even in PC$(T_1,T))$ of cardinality
$\lambda$.  But by the variant EF$^+_{\alpha,\lambda}$-equivalent,
such results are excluded; by it we define our choice
test problem the version of being fat/lean, see Definition \ref{0z.3}.

Why EF$^+$?  See Discussion \ref{1a.8}.

Concerning variants of strongly dependent theories
see  \cite[\S3]{Sh:783},\cite{Sh:863} (and maybe \cite{Sh:F705}), 
most relevant is \cite[\S5]{Sh:863}, part (F).  The best relative 
for us is ``strongly$_4$ dependent", a definition of it is given below
but we delay the treatment to a subsequent paper, \cite{Sh:F918}.
There we also deal with the relevant logics and more.

We prove
here that if $T$ is not strongly stable then $T$ is consistently fat.
More specifically, for every $\mu = \mu^{< \mu} > |T|$ 
there is a $\mu$-complete class forcing
notion $\bbP$ such that in $\bold V^{\bbP}$ the theory $T$ is fat.  
The result holds even for PC$(T_1,T)$.  This gives new cases even for
PC$(T)$ by \ref{0z.5}.

Also if $T$ is unstable or has the DOP or OTOP (see \ref{0z.23}
below or \cite{Sh:c}) \then \, it is fat, i.e. already in $\bold V$.  

Of course, forcing the example is a drawback, but note that for
proving there is no positive theory it is certainly enough.  Hence it
gives us an upper bound on the relevant dividing lines.

On Eherenfeucht-Mostowski models,
see \cite[Ch.III]{Sh:e} or \cite[Ch.VII]{Sh:c} or \cite{Sh:h},
\cite[Ch.III,\S1]{Sh:e}.  I thank a referee for pointing out on earlier
version that Hyttinen-Shelah
\cite{HySh:474} was forgotten hence as Definition \ref{1a.9} was not
yet written, the main result \ref{3c.4} had not said anything new.

I also thank referees for many helpful remarks.

\subsection {Notations and Basic Definitions} \

\begin{definition}
\label{0z.3}
Let $T$ be a complete first order theory.

\noindent
1) We say $T$ is fat \when \, for every
ordinal $\kappa$, for some (regular) cardinality $\lambda > \kappa$ 
there are non-isomorphic models $M_1,M_2$ of $T$
of cardinality $\lambda$ which are 
EF$^+_{\beta,\kappa,\kappa,\lambda}$-equivalent for
every $\beta < \lambda$ (see Definition \ref{1a.9} below).

\noindent
2) If $T$ is not fat, we say it is lean.

\noindent
3) We say the pair $(T,T_1)$ is fat/lean \when \, 
($T_1$ is first order $\supseteq T$ and) PC$(T_1,T) := 
\{M \restriction \tau_T:M$ a model of $T_1\}$ is as above.

\noindent
4) We say $(T,*)$ is fat \when \, for every first order $T_1 \supseteq T$
the pair $(T,T_1)$ is fat.  We say $(T,*)$ is lean otherwise.
\end{definition}

Our claims (mainly \ref{3c.4}) seem to make it clear that some
stable $T$ has NDOT and NOTOP which falls under \ref{3c.4}, but a
referee asks for an example, see \cite[\S5(F)]{Sh:863} for details.

\begin{example}
\label{0z.5} 
1) There is a stable NDOP,NOTOP countable
complete theory which is not strongly dependent; (moreover not 
is not strongly$_4$ stable), see \cite[\S5(G)]{Sh:863}.

\noindent
2) $T = \text{ Th}({}^{\omega_1}(\bbZ_2),E_n)_{n < \omega}$ is as
   above where $\bbZ_2 = \bbZ/2 \bbZ$ as an additive group, $E_n
   = \{(\eta,\nu):\eta,\nu \in {}^{\omega_1}(\bbZ_2)$ are such that
   $\eta \rest (\omega n) = \nu \rest (\omega n)$ where we interpret
   $\bbZ_2$ as the additive group $(\bbZ/2 \bbZ,+,0)$ and
   ${}^{\omega_1}(\bbZ_2)$ as its $\omega_1$-th power as an abelian group.
\end{example}

\centerline {$* \qquad * \qquad *$}

The reader may look at the definitions below only when used.
\begin{definition}
\label{0z.7}
1) Mod$_T(\lambda) = \text{ EC}_T(\lambda)$ is the 
class of models of $T$ of cardinality $\lambda$
and Mod$_T = \text{ EC}_T$ is $\cup\{\text{EC}_T(\lambda):\lambda$ a
cardinality$\}$.

\noindent
2) PC$_\tau(T) = \{M \restriction \tau:M$ a model of $T\}$ where $T$
is a theory or a sentence, in whatever logic, in a vocabulary $\tau_T
\supseteq \tau$; if $\tau = \tau_T$ we may omit $\tau$.

\noindent
3) If $T \subseteq T_1$ are complete first order theories then
 PC$(T_1,T) = \text{ PC}_{\tau(T)}(T_1)$.
\end{definition}

\begin{notation}
\label{0z.11}
1) $\ell g(\bar a)$ is the length of a  sequence $\bar a$.

\noindent
2) $\bar a \trianglelefteq \bar b$ means that $\bar a$ is an initial
   segment of $\bar a$.

\noindent
3) $\bar a \restriction \alpha$ is the unique initial segment of $\bar
   a$ of length $\alpha$ for $\alpha \le \ell g(\bar a)$.
\end{notation}

\begin{definition}
\label{0z.13}
1) For a regular uncountable cardinal $\lambda$ let 
$\check I[\lambda] = \{S \subseteq \lambda$: some
pair $(E,\bar a)$ witnesses $S \in \check I(\lambda)$, see below$\}$. 

\noindent
2) We say that $(E,u)$ is a witness for $S \in \check I[\lambda]$ 
\when \,:
\mn
\begin{enumerate}
\item[$(a)$]   $E$ is a club of the regular cardinal $\lambda$
\sn
\item[$(b)$]   $\bar u = \langle u_\alpha:\alpha < \lambda
\rangle,a_\alpha \subseteq \alpha$ and $\beta \in a_\alpha \Rightarrow
a_\beta = \beta \cap a_\alpha$
\sn
\item[$(c)$]   for every $\delta \in E \cap S,u_\delta$ is an
unbounded subset of $\delta$ of order-type $< \delta$ (and $\delta$ is
a limit ordinal).
\end{enumerate}
\end{definition}

\begin{notation}
\label{0z.17}
1) For a model $M,\bar a \in {}^\alpha M,B \subseteq M$ 
and $\Delta$ a set of formulas, we are interested in
formulas of the form $\varphi(\bar x,\bar y),\bar x = \langle x_i:i <
\alpha\rangle$, so $\alpha$ may be infinite, but the formulas here are
normally first order, so all but finitely many of the $x_i$'s are
dummy variables.

\noindent
1A) $\text{tp}_\Delta(\bar a,B,M) = \{\varphi(\bar x,\bar a):\varphi(\bar
x,\bar y) \in \Delta \text{ and } \bar b \in {}^{\ell g(\bar y)}A 
\text{ and } M \models \varphi[\bar a,\bar b]\}$.

\noindent
2) If $\Delta_{\text{qf}}$ is the set of quantifier-free formulas in
$\bbL(\tau_M)$, we may write tp$_{\text{qf}}$ instead of tp$_\Delta$.

\noindent
3) $\dot I(\lambda,T)$ is the number of isomorphic types of models of
   $T$ of cardinality $\lambda$.

\noindent
4) $\dot I_\tau(\lambda,T)$ is the number of isomorphic types of $M
   \restriction \tau,M$ a model of $T$ of cardinality $\lambda$.

\noindent
5) $\dot I \dot E_\tau(\lambda,T)$ is the supremum of $\{|K|:K
   \subseteq \text{ PC}_\tau(T)$ and $M \in K \Rightarrow \|M\| =
\lambda$ no $M \in K$ has an elementary embeding into any $N \in K 
\backslash \{M\}$, writing $\dot I \dot E_\tau(\lambda,T) =^+ \chi$ we
   mean the supremum is obtained if not said otherwise.

\noindent
6) $\dot I \dot E(\lambda,T) = \dot I \dot E_{\tau(T)}(\lambda,T)$.
\end{notation}

\begin{definition}
\label{0z.23}
Let $T$ be a first order complete theory.

\noindent
1) $T$ has OTOP \when \, $T$ is stable and for some $n,m$ letting
   $\bar x = \langle x_\ell:\ell < n\rangle,\bar y = \langle
   y_\ell:\ell < n\rangle,\bar z = \langle z_\ell:\ell < m\rangle$,
   there are complete types $p(\bar x,\bar y,\bar z)$ such that: for
   every $\lambda$ there is a model $M$ of $T$ and $\bar a_\alpha \in
   {}^n M$ for $\alpha < \lambda$ such that:
\mn
\begin{enumerate}
\item[$(a)$]  $\langle \bar a_\alpha:\alpha < \lambda\rangle$ is an
   indiscernible set
\sn
\item[$(b)$]  for $\alpha \ne B < \lambda$ the type $(p(\bar
a_\alpha,\bar a_b,\bar z)$ is realized in $M$ iff $\alpha < \beta$.
\end{enumerate}
\mn
1A) $T$ has the NOTOP \when \, it is stable but fail the OTOP.

\noindent
2) $T$ has NDOP \when \, $T$ is stable and we can find
$|T|^+$-saturated models $M_\ell$ of $T$ for $\ell \le 3$ such that
   $M_0 \prec M_\ell \prec M_3$ for $\ell=1,2$ and tp$(M_1,M_2)$ does
   not fork over $M_0,M_3$ is $|T|^+$-prime over $M_1 \cup M_2$ but
   not $|T|^+$-minimal over it; equivalently for every $\bar c 
\in {}^{\omega >}(M_3)$
 the type tp$(\bar c,M_1 \cup M_2,M_3)$ is $|T|^+$-isolated but
   there is no infinite $\bold I \subseteq M_3$ which is indiscernible
   over $M_1 \cup M_2$.

\noindent
2A) $T$ has DOP \when \, $T$ is stable and fail to have the NDOP.
\end{definition}

\begin{definition}
\label{0z.25}
1) For a complete first order theory
 $T$, we can say that $\psi$ is a $(\mu,\kappa,T)$-candidate when:
\mn
\begin{enumerate}
\item[$(a)$]   $\psi \in \bbL_{\kappa^+,\omega}(\tau_*)$ for
some vocabulary $\tau_* \supseteq \tau_T$ of cardinality $\le \kappa$
\sn
\item[$(b)$]   PC$_{\tau(T)}(\psi) \subseteq \text{ EC}(T)$
\sn
\item[$(c)$]   for some\footnote{i.e. $\Phi$ proper for
$K^\omega_{\tr}$, i.e. normal trees with $\omega +1$ level, with
linear order on the successor of each node of finite level,
see Definition \ref{2b.1.7}(7) or \cite[Ch.VII]{Sh:c}} $\Phi
\in \Upsilon^{\omega-\tr}_\kappa$ satisfying $\tau_\Phi \supseteq
\tau_\psi$ and EM$({}^{\omega \ge}\lambda,\Phi) \models \psi$ for every
(equivalent some) $\lambda$ and $\Phi$ witness $T$ is not superstable. 
\end{enumerate}
\end{definition}

\bn
Recall that by \cite[Ch.VII]{Sh:c}:
\begin{claim}
\label{0z.29}
If a first order complete theory $T$ is not superstable, \then \, 
for some $\Phi \in \Upsilon^{\omega-{\text{\rm tr}}}_{\tau_2}$, see Definition
\ref{2b.1.7}, $\tau_2 \supseteq \tau(\psi)$ of
cardinality $\kappa,\Phi$ witness $T$ is not superstable, i.e. for
some formulas $\varphi_n(x,\bar y_n) \in \bbL(\tau_T)$, if $I =
{}^\omega \lambda,M = \text{\rm EM}(I,\Phi)$ then for $\eta \in
{}^\omega \lambda,n < \omega$ and $\alpha < \lambda$ we have $M
\models \varphi_n[\bar a_\eta,a_{(\eta \rest n) \char 94 <\alpha>}]$
iff $\alpha = \eta(n)$.
\end{claim}

\begin{definition}
\label{0z.37}
1) For any structure $I$ we say that $\langle \bar a_t:t \in 
I\rangle$ is indiscernible (in the model
$\gC$, over $A$, if $A = \emptyset$ we may omit it) when: 
($\bar a_t \in {}^{\ell g(\bar a_t)}{\gC}$ and)
$\ell g(\bar a_t)$, which is not necessarily finite 
depends only on the quantifier-free type of $t$ in $I$ and:
\mn
\begin{enumerate}
\item[${{}}$]  if $n < \omega$ and $\bar s = \langle
s_0,s_1,\dotsc,s_{n-1}\rangle,\bar t = \langle t_0,\dotsc,t_{n-1}\rangle$
realize the same quantifier-free type in $I$ \then \, $\bar a_{\bar t} := 
\bar a_{t_0} \char 94 \ldots \char 94 \bar a_{t_{n-1}}$ and $\bar a_{\bar
s} = \bar a_{s_0} \char 94 \ldots \char 94 \bar a_{s_{n-1}}$ realizes
the same type (over $A$) in ${\gC}$.
\end{enumerate}
\mn
2) We say that $\langle \bar a_u:u \in [I]^{< \aleph_0}\rangle$ is
indiscernible (in ${\gC}$, over $A$) similarly:
\mn
\begin{enumerate}
\item[${{}}$]  if $n < \omega,w_0,\dotsc,w_{m-1} 
\subseteq \{0,\dotsc,n-1\}$ and $\bar s =
\langle s_\ell:\ell < n\rangle,\bar t = \langle t_\ell:\ell <
n\rangle$ realize the same quantifier-free types in $I$ and $u_\ell =
\{s_k:k \in w_\ell\},v_\ell = \{t_k:k \in w_\ell\}$ \then \, $\bar a_{u_0}
\char 94 \ldots \char 94 \bar a_{u_{n-1}},\bar a_{v_0} \char 94 \ldots
\char 94 \bar a_{v_{n-1}}$ realize the same type in ${\frak C}$ (over $A$).
\end{enumerate}
\mn
3) If $I$ is a linear order then we let 
incr$({}^{\alpha} I) = \text{ incr}_\alpha(I) = 
\text{ incr}(\alpha,I)$ be 
$\{\rho:\rho$ is an increasing sequence of length $\alpha$ of members 
of $I\}$; similarly incr$({}^{\alpha >}I)$-incr$_{< \alpha}(I) =
\text{ incr}(< \alpha,I) := \cup\{\text{incr}_\beta(I):\beta < \alpha\}$.  So 
instead $[I]^{< \aleph_0}$ we may use incr$_{<\omega}(I)$; clearly the
difference is notational only.
\end{definition}
\newpage

\section {Games, equivalences and questions} 

What is the meaning in using EF$^+_{\alpha,\lambda}$?  Consider for
various $\gamma$'s the game $\Game_{\gamma,\lambda}(M_1,M_2)$ where
$M_1,M_2 \in \text{ Mod}_T(\lambda),T$ complete first order 
$\bbL(\tau)$-theory.  During a play we can consider dependence relations
on ``short" sequences from $M_\ell$ (where $\le 2^{|\tau| + \aleph_0}$ is
the default value), definable in a suitable sense.  So if $T$ is a
well understood unsuperstable $T$ like Th$({}^\omega \omega,E_n)_{n <
\omega}$ with $E_n := \{(\eta,\nu):\eta,\nu \in {}^\omega \omega$ and
$\eta \restriction n = \nu \restriction n\}$, \then \, even for 
$\gamma = \omega +2$ we have $E^+_{\gamma,\lambda}$- equivalence
implies being isomorphic.  This fits the thesis:

\begin{thesis}
\label{1a.2}
The desirable dichotomy characterized, on the
family of first order $T$, by the property ``$M_1,M_2 \in 
\text{ Mod}_T(\lambda)$ are long game EF-equivalent iff they are isomorphic",
is quite similar to the one in \cite[Ch.XIII]{Sh:c}; the structure side is
e.g.: $T$ is stable and every $M \in \text{ Mod}_T$ is prime over some
$\cup\{M_\eta:\eta \in I\}$, where ${\cT}$ is a subtree of
${}^{\kappa_r(T) >}\|M\|$ and $\eta \triangleleft \nu \Rightarrow M_\eta
\prec M_\nu \prec M,\|M_\eta\| \le 2^{|T|}$ and $\eta \triangleleft \nu
\in \cT \Rightarrow \text{ tp}(M_\nu,\cup\{M_\rho:\rho \in \cT,\rho
\restriction (\ell g(\nu)+1) \ne \eta \restriction (\ell g(\nu)+1))$
does not fork over $M_\nu$, i.e. $\bar M = \langle M_\eta:\eta \in
{\cT}\rangle$ is a non-forking tree of models with $\le
\kappa_r(T)$ many levels.  

We think the right (variant of the) question is from \ref{1a.4}.  Probably a
reasonable analog is the situation in \cite[Ch.XII,XIII]{Sh:c}: the
original question was on the function $\lambda \mapsto \dot
I(\lambda,T)$, the number of non-isomorphic models; but the answer is
more transparent for $\lambda \mapsto \dot I \dot E(\lambda,T)$.

If $\lambda = \mu^+,\mu = \mu^{|T|} = \text{ cf}(\mu),T = \text{
Th}({}^\omega \omega,E_n)_{n < \omega}$ then by Hyttinen-Shelah
\cite[Th4.4]{HySh:474}; for $\gamma \ge \mu \omega$ we get equivalence
$\Rightarrow$ isomorphic, but not for $\gamma < \mu \omega$; now
\ref{1a.19} is parallel to that.  This seems to
indicate that EF$^+_{\gamma,\lambda}$ is suitable for the questions we
are asking: it uses the game EF$^+$, which is more complicated but the
length of the game is much ``smaller" in the relevant results.
\end{thesis}

So the natural question concerning such equivalences is (see
\cite{Sh:c}, \cite{Sh:E53}):

\begin{question}
\label{1a.4}
Classify first order complete $T$, or at least the countable ones by:
\medskip

\noindent
\underline{Version (A)$_1$}:  For every ordinal $\alpha$, there are a
cardinal $\lambda$ and 
non-isomorphic $M_1,M_2 \in \text{ Mod}_T(\lambda)$ which
are EF$^+_{\alpha,\lambda}$-equivalent (at least, e.g. in some $\bold
V^{\bbP},\bbP$ is $(2^{|T|+|\alpha|})^+$-complete forcing notion).
\medskip

\noindent
\underline{Version (A)$_0$}:  Similar version for EF$_{\alpha,\lambda}$.
\medskip

\noindent
\underline{Version (B)$_1$}:  For every 
cardinal $\kappa > |T|$ and vocabulary $\tau_1
\supseteq \tau_T$ and $\psi \in \bbL_{\kappa,\omega}(\tau_1)$ such
that PC$_\tau(\psi) \subseteq \text{ EC}_T$ has members of arbitrarily
large cardinality we have $(a) \Rightarrow (b)$ where
\mn
\begin{enumerate}
\item[$(a)$]   for every cardinal $\mu$ in PC$_\tau(\psi) := \{M
\restriction \tau:M$ a model of $\psi\}$ there is a $\mu$-saturated member
\sn
\item[$(b)$]   for every $\alpha$ for arbitrarily large 
$\lambda$ there are $M_1,M_2 \in
\text{ PC}_\tau(\psi)$ of cardinality $\lambda$ with 
non-isomorphic $\tau$-reducts which are EF$^+_{\alpha,\lambda}$-equivalent.
\end{enumerate}
\medskip

\noindent
\underline{Version (B)$_0$}:  Like (B)$_1$ for EF$_{\alpha,\lambda}$.
\medskip

\noindent
\underline{Version (C)$_1$}:  Like (B)$_1$ using $\psi = \wedge T_1$ where
$T_1$ is first order $\supseteq T$.
\medskip

\noindent
\underline{Version (C)$_0$}:  Like (B)$_0$ using $\psi = \wedge T_1$ where
$T_1$ is a first order $\supseteq T$.
\end{question}

\begin{discussion}
\label{1a.5}
1) For reasons to prefer version (B) over (C) - see \cite{Sh:E53}.

\noindent
2) Now by the works quoted above, (see \cite[3.19]{HySh:529} quoted in
\ref{4d.2} below):  $T$ satisfies (A)$_0$ iff $T$ is
superstable NDOP, OTOP iff (B)$_0$.  Of course if we change the order
of the quantifier (to ``for aribitrarily large some 
$\lambda$ for every $\alpha < \lambda$,...") this is not 
so, but we believe solving (A)$_1$ and/or
(B)$_1$ will eventually do much also for this.
\end{discussion}

So all this means
\begin{conjecture}
\label{1a.6}
1) For a complete (first order) $T$ the following are equivalent:
\mn
\begin{enumerate}
\item[$(a)$]   for every ordinal $\alpha$ for some $\lambda$ there
are non-isomorphic, EF$^+_{\alpha,\lambda}$-equivalent models $M_1,M_2
\in \text{ EC}_T(N)$
\sn
\item[$(b)$]   for arbitrarily large $\lambda$ for every $\alpha <
\lambda$ there are non-isomorphisms, EF$^+_{\alpha,\lambda}$-equivalent
models $M_1,M_2 \in \text{ EC}_T(\lambda)$
\sn 
\item[$(c)$]  for every large enough regular $\lambda$ there are 
non-isomorphisms $M_1,M_2 \in \text{ EC}_T(\lambda)$ which are 
$E_{\alpha,\lambda}$-equivalent for every $\alpha < \lambda$.
\end{enumerate}
\mn
2) Similarly for ``some $T_1 \supseteq T$, PC$(T_1,T)$ is lean.

We conjecture that proving that if we prove that a (countable) fat $T$
is close enough to superstable, will enable us to generalize proofs in
\cite[Ch.XII]{Sh:c} only now the tree has $\le \omega_1$ levels rather
than $\omega$.

We can also return to the ordinals $\alpha \in (\lambda \omega,\lambda^+)$.
\end{conjecture}

\centerline {$* \qquad * \qquad *$}

Now we shall actually look at the games.
\begin{definition}
\label{1a.7}
1) We say that $M_1,M_2$ are EF$_\alpha$-equivalent 
if $M_1,M_2$ are models (with same vocabulary) and $\alpha$ is an ordinal
such that the isomorphism player has a winning strategy in the game
${\cG}^\alpha_1(M_1,M_2)$ defined below.

\noindent
1A) Replacing $\alpha$ by $< \alpha$ means: for every $\beta <
\alpha$; similarly below.

\noindent
2) We say that $M_1,M_2$ are EF$_{\alpha,\mu}$-equivalent or
${\cG}^\alpha_\mu$-equivalent \when \,
$M_1,M_2$ are models with the same vocabulary, $\alpha$ an ordinal,
$\mu$ a cardinal such that the isomorphism
player has a winning strategy in the game ${\cG}^\alpha_\mu(M_1,M_2)$
defined below.

\noindent
3) For $M_1,M_2,\alpha,\mu$ as above and partial isomorphism
 $f$ from $M_1$ into $M_2$ we define the game
 ${\cG}^\alpha_\mu(f,M_1,M_2)$ between the players ISO, the
 isomorphism player and AIS, the anti-isomorphism player as follows:
\mn
\begin{enumerate}
\item[$(a)$]   a play lasts $\alpha$ moves
\sn
\item[$(b)$]  after $\beta$ moves a partial isomorphism $f_\beta$
 from $M_1$ into $M_2$ is chosen, increasing continuous with $\beta$
\sn
\item[$(c)$]  in the $(\beta +1)$-th move, the player AIS chooses $A_{\beta,1}
\subseteq M_1,A_{\beta,2} \subseteq M_2$ such that $|A_{\beta,1}| +
|A_{\beta,2}| < 1 + \mu$ and then the player ISO chooses $f_{\beta +1}
\supseteq f_\beta$ such that 
$A_{\beta,1} \subseteq \text{ Dom}(f_{\beta +1})$ and $A_{\beta,2} \subseteq$
Rang$(f_{\beta +1})$
\sn
\item[$(d)$]   if $\beta=0$, ISO chooses $f_0=f$; if $\beta$ is a
limit ordinal ISO chooses $f_\beta = \cup\{f_\gamma:\gamma < \beta\}$.
\end{enumerate}
\mn
The ISO player loses if he had no legal move for some $\beta <
\alpha$, otherwise he wins the play. 

\noindent
4) If $f = \emptyset$ we may write ${\cG}^\alpha_\mu(M_1,M_2)$.  If $\mu$
is 1 we may omit it.  We may write $\le \mu$ instead of $\mu^+$.
\end{definition}

\begin{discussion}
\label{1a.8}
1) Why do we need EF$^+$?  

First, if we like a parallel of \cite[Ch.XIII]{Sh:c}, i.e. a game in
which set of small cardinality are chosen, say $|T|$ or $2^{|T|}$ or
whatever rather than just $< \lambda = \|M_\ell\|$, clearly
EF$_{\alpha,\mu}$ cannot help.

\noindent
2) Also, consider $\lambda = \mu^+,\mu = \text{ cf}(\mu) > |T|$ and an
   ordinal $\alpha < \lambda$ and ask for which $T$: for any two
   models $M_1,M_2$ of $T$ of cardinality $\lambda$,
   EF$_{\alpha,\lambda}$-equivalence implies isomorphisms?  (The
   EF$_{\alpha,\lambda}$-equivalence means that the isomorphism player
   wins in the game of length $\alpha$, in each step adding $\le \mu$
   elements to the domain and range of the partial isomorphism.)

Now we know (by earlier works, see \ref{4d.23}) for countable $T$
that if $\alpha \in [\omega,\mu \times \omega]$ that the answer (for
the pair $(\alpha,\lambda)$) is as in the main gap for $\dot I \dot E$ 
($T$ superstable with
NDOP and NOTOP).  But for larger $\alpha < \lambda$ this is not so, as
e.g. for the prototypical stable unsuperstable $T$ for $\alpha = \mu
\times (\omega +2)$ we get yes, ``it is low".

\noindent
3) Looking at the reason for this, i.e. why we need $\mu \times (\omega +2)$
 moves, not $(\omega +2)$ moves we formulate EF$^+$.  We think that
   with EF$^+_{\alpha,\theta,\mu,\lambda}$ for small
   $\alpha,\theta,\mu$ and just $\lambda = \|M_\ell\|$ we get the
   desired dichotomy.  In general, we expect the results will be
   robust under choosing such an exact game; and will resolve the case
   $\alpha \in (\mu \times (\omega,2),\lambda)$ case above.

\noindent
4) More specifically, the reason EF$_{\alpha,\lambda}$-equivalence
   does not imply isomorphisms for $M_1,M_2 \in 
\text{ EC}_\lambda(T)$, even in the case $T = 
\text{ Th}({}^\omega \omega,E_n)_{m < \omega}$, is that: 
assume we fix a winning strategy {\bf st} for
   ${\cG}_{\alpha,\lambda}(M_1,M_2)$, if we let $\langle
   a^\ell_\alpha/E^{M_\ell}_1:\alpha < \lambda\rangle$ list
   $M_\ell/E^{M_\ell}_1$ and $\bold R = \{(\alpha,\beta)$: in 
  some short initial segment $\bold x$ of a play of 
${\cG}_{\alpha,\lambda}(M_1,M_2)$ in which the player ISO uses the
strategy {\bf st}, we have $f^{\bold x}_\alpha(a^1_\alpha)
E^{M_2}_1 a^2_\beta\}$, we have to find a
   function $h$ from $\lambda$ onto $\lambda$ whose graph is
   $\subseteq \bold R$.

Now being in a winning position is enough to show the existence of
such $h$, only when the game is long enough.  For EF$^+_{\alpha,\theta}$
this is different.

\noindent
5) Note: we use the case $k=1$ from \ref{1a.9}.  If we shall have
good structure theorems then even $k=2$ is O.K.  For $k = \kappa$ 
it expresses the
logic in \cite[Ch.XIII]{Sh:c} when we add the game quantifier of
appropriate length.

\noindent
6) Of course, the case $k=0$ is easier for ISO then the case $k=2$
   which is easier than $k=1$, so the relevant implications holds.
\end{discussion}

\begin{definition}
\label{1a.9}
1) For $k \in \{0,1,2\}$ the models $M_1,M_2$ are 
EF$^{+,k}_{\gamma,\theta,\mu,\lambda}$-equivalent, but if
$k=1$ we may omit it, \when \, the isomorphism player, ISO,
has a winning strategy in the game
${\cG}^k_{\gamma,\theta,\mu,\lambda}(M_1,M_2)$ defined below.

We always assume $\aleph_0 \le \theta \le \mu$.
If $\mu = \text{ min}\{\|M_1\|,\|M_2\|\}$ then we may omit it.  
If also $\theta = (2^{|\tau(M_\ell)|+\aleph_0})^+$ we may omit it, too.

\noindent
2) For an ordinal $\gamma$,
cardinals $\theta \le \mu$, vocabulary $\tau$ and
$\tau$-models $M_1,M_2$ and partial isomorphism $f$ from $M_1$ to $M_2$,
we define a game ${\cG}^k = 
{\cG}^{+,k}_{\gamma,\theta,\mu,\lambda}(f,M_1,M_2)$,
between the player ISO
(isomorphism) and AIS (anti-isomorphism).  

A play last $\gamma$ moves; in the $\beta$-th move a partial
isomorphism $f_\beta$ from $M_1$ to $M_2$ is chosen by ISO, extending
$f_\alpha$ for $\alpha < \beta$ such that $f_0 = f$ and for limit
$\beta$ we have $f_\beta = \cup\{f_\alpha:\alpha < \beta\}$ and for
every $\beta < \alpha$ the set Dom$(f_{\beta +1}) \backslash
\text{ Dom}(f_\beta)$ has cardinality $< 1 + \mu$; let $f^\ell_\beta$
be $f_\beta$ if $\ell = 1,f^{-1}_\beta$ if $\ell=2$.  

During a play, the player ISO loses if he has no legal move and he
wins in the end of the play iff he always had a legal move.

In the $(\beta +1)$-th move, the AIS player does one of the following cases:
\bigskip

\noindent
\underline{Case 1}:  The AIS player chooses $A_\ell = A^\ell_\beta 
\subseteq M_\ell$ for $\ell=1,2$ such that $|A_1| + |A_2| < 1 +\mu$ and 
then ISO chooses $f_\beta$ as above such that $A_\ell 
\subseteq \text{ Dom}(f^\ell_\beta)$ for $\ell = 1,2$.
\bigskip

\noindent
\underline{Case 2}:  First the AIS player 
chooses\footnote{note that for $k=0,1$ we require ``$\bbL(\tau_T)$-definable
$\bold R_\ell$ such that $f$ maps the definition of $\bold R_1$ to the
one of $\bold R_2$"; moreover we expect that we can demand it
is as in the case of using regular types.} a
pre-dependence relation $\bold R_\ell$ on ${}^{\theta >}(M_\ell)$ 
(see Definition \ref{1a.15} below) and ${\cA}_\ell 
\subseteq {}^\varepsilon (M_\ell)$ of cardinality $\le \lambda$ for 
$\ell=1,2$ such that:
\mn
\begin{enumerate}
\item[$\odot$]   $(a) \quad$ if $k=0$ then $\bold R_\ell =
[{}^{\theta >}(M_\ell)]^{< \aleph_0}$, so really an empty case
\sn
\item[${{}}$]   $(b) \quad$ if $k=1,2$ then $\bold R_\ell$ is a
1-dependence relation (see \ref{1a.15}(4)(b)(B) 

\hskip25pt below)
\sn
\item[${{}}$]  $(c) \quad$ if $k=1,2$ and $\ell=1,2$ and $n <
\omega$ and $\bar a_0,\dotsc,\bar a_{n-1} \in {}^\varepsilon(M_\ell)$
then

\hskip25pt  the truth value of $\{\bar a_0,\dotsc,\bar a_{n-1}\} 
\in \bold R_\ell$ depends just on the 

\hskip25pt  complete first order type which $\langle \bar a_0,\dotsc,\bar
a_{n-1}\rangle$ realizes on 

\hskip25pt Dom$(f^\ell_\beta)$ inside the model $M_\ell$.
\end{enumerate}
\mn
Second, the ISO does one of the following:
\bigskip

\noindent
\underline{Subcase 2A}:  First, assume $k=2$.  The player ISO chooses 
$\langle(\bar a^1_\zeta,\bar a^2_\zeta):\zeta < \lambda\rangle$ such
that for $\ell=1,2$:
\mn
\begin{enumerate}
\item[$(\alpha)$]   for each $\zeta < \lambda$ for some $\varepsilon <
\theta$ we have $\bar a^\ell_\zeta \in {}^{\varepsilon}(M_\ell)$ 
\sn
\item[$(\beta)$]   $\langle \bar a^\ell_\zeta:\zeta <
\lambda\rangle$ is independent for $\bold R_\ell$
\sn
\item[$(\gamma)$]    each $\bar a \in {\cA}_\ell$ does 
$\bold R_\ell$-depend on $\{\bar a^\ell_\zeta:\zeta < \lambda\}$. 
\end{enumerate}
\mn
Then AIS chooses $\zeta < \lambda$ and ISO chooses $f_{\beta +1} \supseteq
f_\beta$ such that $f_\beta(\bar a^1_\zeta) = \bar a^2_\zeta$.

Second, assume $k=1$.  Then the ISO player chooses
equivalence relations $E_\ell$ on ${}^{\theta >}(M_\ell)$ which the
dependence relation, i.e. $E_{\bold R_\ell}$, i.e. \ref{1a.15}(6)
 and equality of length refine for $\ell=1,2$ and
choose a function $h$ from the family of $E_1$-equivalence classes
onto the family of $E_2$-equivalence classes which preserve
cardinality up to $\lambda$; that is, if $h(\bar a_1/E_1) = \bar
a_2/E_2$ then $\ell g(\bar a_1) = \ell g(\bar a_2) \and \text{
min}\{\text{dim}(\bar a_1/E_1),\lambda) = \text{ min}\{\text{dim}(\bar
a_2/E_2),\lambda\}$.

Then the AIS player chooses a pair $(\bar a_1,\bar a_2)$ such that
$\bar a_\ell \in {}^{\theta >}(M_\ell)$ for $\ell=1,2$ such that
$h(\bar a_1/E_1) = (\bar a_2/E_2)$ and ISO has to choose 
$f_{\beta +1} \supseteq f_\beta$ such that $f(\bar a_1) = \bar a_2$.  
\bigskip

\noindent
\underline{Subcase 2B}:  The 
player ISO chooses $f_{\beta +1} \supseteq f_\beta$
as required such that for some $n < \omega$ and $\bar a^1_\ell \in
{}^\varepsilon\text{Dom}(f_\beta)$ for $\ell \le n$ we have: 
$\{\bar a^1_0,\dotsc,\bar a^1_{n-1}\}$ is $\bold R_1$-dependent 
\Iff \, $\{f_\beta(\bar a^1_0),\dotsc,
f_\beta(\bar a^1_{n-1})\}$ is not $\bold R_2$-dependent.
\end{definition}

\begin{definition}
\label{1a.15}
1) We say $\bold R$ is a pre-dependence relation 
on $X$ \when \, $\bold R$ is a subset of $[X]^{< \aleph_0}$.

\noindent
2) For $X,\bold R$ as above, we say $Y \subseteq X$ is $\bold
R$-independent when $[Y]^{< \aleph_0} \cap \bold R = \emptyset$; of
course, an index set with repetitions is considered dependent.

\noindent
3) We say $\bold R$ or $(X,\bold R)$ has character $\le \kappa$ \when \,
for every $\bold R$-independent $Y \subseteq X$ and $\{x\} \subseteq
X$ for some $Z \in [Y]^{< \kappa}$ the set $(Y \backslash Z) \cup
\{x\}$ is $\bold R$-independent.

\noindent
4) We say that $\bold R$ is a $k$-dependence relation on $X$ (if $k=1$
we may omit it) when:
\mn
\begin{enumerate}
\item[$(a)$]   $\bold R$ is a subset of $[X]^{< \aleph_0}$
\sn
\item[$(b)$]   $(\alpha) \quad$ if $k=0$ then $\bold R = [X]^{< \aleph_0}$
\sn
\item[${{}}$]  $(\beta) \quad$ if $k=1$ then $\bold R$-independence
satisfies the exchange principle

\hskip25pt  (so dimension is well defined, as for regular types).
\end{enumerate}
\mn
5) We say $R$ is trivial \when \, for every $Y \subseteq X,Y$ is
$\bold R$-independent iff every $Z \subseteq [Y]^{\le 2}$, is
   $\bold R$-independent.

\noindent
6) For $\bold R$ as in $(a),(b)(\beta)$ let $E_{\bold R} =
   \{\{x_1,x_2\}:x_1 =x_2 \in X$ on $\{\lambda_1\},\{x_2\} \in \bold
   R$ or $\{x_1,x_2\} \in \bold R \wedge \{x_1\} \notin \bold R \wedge
   \{x_2\} \notin \bold R\}$ is an equivalence relation on $X$;
   pedantically we should write $E_X,\bold R$.  
\end{definition}

\begin{claim}
\label{1a.19}
$M_1,M_2$ are isomorphic \when \,:
\mn
\begin{enumerate}
\item[$(a)$]  $M_1,M_2$ are models of $T$ of cardinality $\lambda$
\sn
\item[$(b)$]   $M_1,M_2$ are 
{\rm EF}$^+_{\omega +2,\aleph_0,\aleph_0,\lambda}$-equivalent
\sn
\item[$(c)$]   $T = \text{\rm Th}({}^\omega \omega,E_n)_{n <
\omega}$ and $E_n = \{(\eta,\nu):\eta \in {}^\omega \omega,\nu \in
{}^\omega \omega$ and $\eta \rest n = \nu \rest n\}$.
\end{enumerate}
\end{claim}

\begin{PROOF}{\ref{1a.19}}
\bigskip

\noindent
\underline{Step A}:  We choose a winning strategy {\bf st} of
the isomorphism player in the game
${\cG}_{\omega +2,\aleph_0,\aleph_0,\lambda}(M_1,M_2)$.
\bigskip

\noindent
\underline{Step B}:  By the choice of $T$ for $\ell=1,2$ we can find
${\cT}_\ell,\bar{\bold a}_\ell$ such that:
\mn
\begin{enumerate}
\item[$\boxtimes_\ell$]  $(a) \quad {\cT}_\ell$ is a subtree
of ${}^{\omega >}\lambda$
\sn
\item[${{}}$]  $(b) \quad \bar{\bold a}_\ell = \langle
a^\ell_\eta:\eta \in {\cT}_\ell\rangle$
\sn
\item[${{}}$]  $(c) \quad a^\ell_\eta \in M_\ell$
\sn
\item[${{}}$]   $(d) \quad$ if $\eta \in {\cT}_\ell$ and $\ell
g(\eta) = n$ then $\langle a^\ell_\nu/E^{M_\ell}_{n+1}:\nu \in \text{
suc}_{{\cT}_\ell}(\eta)\rangle$ list

\hskip25pt  $\{b/E^{M_\ell}_{n+1}:b \in M_\ell,b \in 
a^\ell_\eta/E^{M_\ell}_n\}$ without repetitions.
\end{enumerate}
\mn
Let ${\cT}_{\ell,n} = \{\eta \in {\cT}_\ell:\ell g(\eta) = n\}$
and let ${\cT}_{\ell,\omega} = \{\eta \in {}^\omega \lambda:\eta
\restriction n \in {\cT}_\ell$ for every $n < \omega\}$.

Lastly, let $\bar \mu_\ell = \langle \mu^\ell_\eta:\eta \in 
{\cT}_{\ell,\omega}\rangle$, where

\[
\mu^\ell_\eta = |\{b \in M_\ell:b \in a^\ell_{\eta \restriction
n}/E^{M_\ell} \text{ for every } n < \omega\}|.
\]
\bigskip

\noindent
\underline{Step C}:

Clearly
\mn
\begin{enumerate}
\item[$\boxplus$]   $M_1,M_2$ are isomorphic \Iff \, there is an
isomorphism $h$ from ${\cT}_1$ onto ${\cT}_2$ (i.e. $h$ maps
${\cT}_{1,n}$ onto ${\cT}_{2,n},h$ preserves the length, 
$\eta \triangleleft
\nu$ and $\eta \ntriangleleft \nu$) such that letting $h_n = h
\restriction {\cT}_{1,n}$ and $h_\omega$ be the mapping from 
${\cT}_{1,\omega}$ onto ${\cT}_{2,\omega}$ which $h$ induces (so
$h_\omega(\eta) = \bigcup\limits_{n < \omega} 
h_n(\eta \restriction n))$ we have
$\eta \in {\cT}_{1,\omega} \Rightarrow \mu^1_\eta =
\mu^2_{h_\omega(\eta)}$. 
\end{enumerate}
\bigskip

\noindent
\underline{Step D}:

By induction on $n$ we choose $h_n,\bar{\bold x}_n$ such that
\mn
\begin{enumerate}
\item[$\circledast$]  $(a) \quad h_n$ is a one-to-one mapping from 
${\cT}_{1,n}$ onto ${\cT}_{2,n}$
\sn
\item[${{}}$]   $(b) \quad$ if $m < n$ and $\eta \in {\cT}_{1,n}$
then $h_m(\eta \restriction m) = (h_n(\eta)) \restriction m$
\sn
\item[${{}}$]  $(c) \quad \bar{\bold x}_n = \langle \bold
x^n_\eta:\eta \in {\cT}_{1,n}\rangle$
\sn
\item[${{}}$]   $(d) \quad (\alpha) \quad \bold x^n_\eta$ is an
initial segment of a play of the game
$\Game_{\omega +2,\aleph_0,\aleph_0,\lambda}(M_1,M_2)$
\sn
\item[${{}}$]   $\quad \quad (\beta) \quad$ in $\bold x^n_\eta$
only finitely many moves have been played

\hskip33pt (can specify), the last one is $m(\bold x^n_\eta)$
\sn
\item[${{}}$]   $\quad \quad (\gamma) \quad$ in $\bold x^n_\eta$, 
the player ISO uses his winning strategy {\bf st}
\sn
\item[${{}}$]   $(e) \quad$ if $\eta_1 \in {\cT}_{1,n}$ and $\eta_2 =
h_n(\eta_1)$, \then \, for some $b_1 \in 
\text{ Dom}(f^{\bold x^n_\eta}_{m(\bold x^n_\eta)})$

\hskip25pt  we have
\sn
\item[${{}}$]  $\quad \quad (\alpha) \quad b_1  \in a^1_\eta/E^{M_1}_n$
\sn
\item[${{}}$]   $\quad \quad (\beta) \quad 
f^{\bold x^n_\eta}_{m(\bold x^n_\eta)}(b_1) \in a^2_{h_n(\eta)}/E^{M_2}_n$
\sn
\item[${{}}$]   $(f) \quad$ if $\nu \triangleleft \eta \in 
{\cT}_{1,n}$ \then \, $\bold x^{\ell g(\nu)}_\nu$ is an initial segment
of $\bold x^n_\eta$.
\end{enumerate}
\mn
Why can we carry the induction?
\bigskip

\noindent
\underline{For $n=0$}:

Note that $h_0$ is uniquely determined.  As for $\bold x^0_{<>}$, any
$\bold x$ as in $\circledast(d)$ is O.K., as long as at least one move
was done (note that $E^{M_\ell}_0$ has one and only one equivalence
class.
\bigskip

\noindent
\underline{For $n=m+1$}:  So $h_m,\bar{\bold x}_m$ has been chosen.

Let $\eta_1 \in {\cT}_{1,m}$ and let $\eta_2 = h_m(\eta_1)$ and

\begin{equation*}
\begin{array}{clcr}
\bold F_{\eta_1} := \{(\nu_1,\nu_2):&\nu_1 \in 
\text{ suc}_{{\cT}_1}(\eta_1),\nu_2 \in \text{ suc}_{{\cT}_2}(\eta_2) \\
  &\text{ and there is } \bold x \text{ as in } \circledast(d) \text{
such that} \\
  &\bold x^m_{\eta_1} \text{ is an initial segment of } \bold x
\text{ and for some} \\
  &b_1 \in \text{ Dom}(f^{\bold x}_{m(\bold x)}) \text{ we have} \\
  &b_1 \in a^1_{\nu_1}/E^{M_1}_n \text{ and } f^{\bold x}(b_1) \in
a^2_{\nu_2}/E^{M_2}_n\}.
\end{array}
\end{equation*}

\mn
Now
\mn
\begin{enumerate}
\item[$\odot$]   to do the induction step, it suffices to prove
that: if $\eta_1 \in {\cT}_{1,m}$ then there is a one-to-one function
$h_{n,\eta_1}$ from suc$_{{\cT}_1}(\eta_1)$ onto 
suc$_{{\cT}_2}(\eta_2)$ such that $\nu \in \text{ suc}_{{\cT}_1}(\eta_1)
\Rightarrow (\nu,h_{n,\eta_1}(\nu)) \in \bold F_{\eta_1}$.
\end{enumerate}
\mn
However by case 2 in Definition \ref{1a.9} this holds.
\bigskip

\noindent
\underline{Stage E}:  

So we can find $\langle h_n:n < \omega\rangle,\langle \bold
x_\eta:\eta \in {\cT}_1\rangle$ as in $\circledast$.  Let $h :=
\cup\{h_n:n < \omega\}$, clearly it is an isomorphism from 
${\cT}_1$ onto ${\cT}_2$ and $h_\omega$ is well defined, see
$\boxplus$ from Stage C.

So it is enough to check the sufficient condition for $M_1 \cong M_2$
then, i.e. $\eta \in {\cT}_{1,\omega} \Rightarrow \mu_{1,\eta} =
\mu_{2,h_\omega(\eta)}$.  But if $\eta \in {\cT}_{1,\omega}$ then
$\langle \bold x_{\eta \restriction n}:n < \omega\rangle$ is a
sequence of initial segments of a play of ${\cG}$ with ISO using his
winning strategy {\bf st}, increasing with $n$, each with 
finitely many moves.  So $\bold x_\eta$, defined as
the limit $\langle \bold x_{\eta \restriction n}:n <\omega\rangle$, 
is an initial segment of the play ${\cG}$, with $\le
\omega$ moves and $f^{\bold x_\eta}_{\bold m(\bold x_\eta)} = 
\cup\{f^{\bold x_{\eta \rest n}}_{m(\bold x_{\eta \rest n)}}:n <\omega\}$.

Clearly $n < \omega \Rightarrow f(a^1_{\eta \restriction n})E^{M_2}_n
a^2_{h_n(\eta \restriction n)}$.  As we have one move left and can use case 2
in Definition \ref{1a.9}(2) we are done. 
\end{PROOF}

\noindent
The following claim says that the game in \ref{1a.7}, \ref{1a.9}
when $\lambda = \mu^+,\alpha < \lambda$ divisible enough are
equivalent, i.e. the ISO player wins one iff he wins the other.
\begin{claim}
\label{1a.23}
1) $M_1,M_2$ are {\rm EF}$^+_{\gamma,\theta,\mu,\lambda}$-equivalent \when \,:
\mn
\begin{enumerate}
\item[$(a)$]   $M_1,M_2$ are $\tau$-models
\sn
\item[$(b)$]  $\lambda = \lambda^+_1,\lambda_1 \ge \mu$ and $\theta
\le \mu \le \lambda$ and $\gamma \le \mu$ and {\rm cf}$(\mu) < \mu
\Rightarrow \lambda_1 > \mu$ and $\lambda \in \check I[\lambda]$, see
Definition \ref{0z.13}
\sn
\item[$(c)$]   $M_1,M_2$ are {\rm EF}$_{\gamma(*),\mu}$-equivalent
where $\gamma(*) = \lambda_1 \times \gamma$ (see Definition
\ref{1a.7}(2))
\sn
\item[$(d)$]   $\|M_\ell\| = \lambda = \lambda^{< \theta}$ for $\ell=1,2$
\end{enumerate}
\mn
2) $M_1,M_2$ are {\rm EF}$^+_{\gamma,\mu}$-equivalent \when \, they are
{\rm EF}$_{\gamma,\theta,\mu,\lambda}$-equivalent.

\noindent
3) $M_1,M_2$ are {\rm EF}$^+_{\gamma_1,\theta_1,\mu_1,\lambda_1}$-equivalent 
\when \, they are {\rm EF}$^+_{\gamma_2,\theta_2,\mu_2,\lambda_2}$-equivalent
and $\gamma_1 \le \gamma_2,\theta_1 \le \theta_2,\mu_1 \le
\mu_2,\lambda_1 \le \lambda_2$.
\end{claim}

\begin{PROOF}{\ref{1a.22}}
1) First, we do not save on $\gamma(*)$, say use $\lambda_1
\times \lambda _1 \times \gamma$.

Let {\bf st} be a winning strategy of the ISO player in the game
${\cG}^{\gamma(*)}_\mu$.  We try to use
it as a winning strategy of the ISO player in the game
${\cG}_{\gamma,\theta,\mu,\lambda}(M_1,M_0)$.  
Well, the $f^{\bold x}_\alpha$ may have too large a 
domain, so ``on the side" in the $\beta$-th move ISO
play $\bold x_\beta$ for ${\cG}_{\gamma,\theta,\mu,\lambda}$ and
$A^1_\beta \subseteq \text{ Dom}(f^{\bold x_\beta})$ of cardinality $<
\mu$ (or $\le \mu$ if $\mu > \text{ cf}(\mu) \wedge \beta \ge \text{
cf}(\mu))$ and he actually plays $f^{\bold x_\beta} \restriction
A^1_\beta$, i.e. is an initial segment of a play of 
${\cG}_{\gamma,\mu}$ of length $\beta$ in which the ISO player uses the
strategy {\bf st} such that $[\beta_1 < \beta \Rightarrow \bold
x_{\beta_1}$ is an initial segment of $\bold x_\beta]$.

The only problem is when $\beta = \alpha +1$ and in Definition
\ref{1a.9}, Case 2 occurs, i.e. with the AIS player choosing
$\bold R^1_\beta,\bold R^2_\beta$.  We may for notational simplicity choose 
$\varepsilon < \theta$ and deal only with $A_\ell \cap {}^\varepsilon(M_\ell)$
for $\ell=1,2$.

We can consider $\bold x_\beta$ extending $\bold x_\alpha$; if it is as
required in subcase (2B) of Definition \ref{1a.9} we are done.  Let

\begin{equation*}
\begin{array}{clcr}
\bold F^1_\beta = \{(\bar a_1,\bar a_2):&\text{ for some } \varepsilon
< \theta,\bar a_\ell \in {}^\varepsilon(M_\ell) \text{ for } \ell =1,2 \\
  &\text{ and there is a candidate } \bold x_\beta \text{ for the} \\
  &\beta\text{-th move such that } f^{\bold x_\beta}(\bar a_1) = \bar a_2\}.
\end{array}
\end{equation*}

\mn
Let

\[
\bold F^2_\beta = \{(\bar a_2,\bar a_1):(\bar a_1,\bar a_2) \in \bold
F^1_\beta\}
\]

\[
{\cA}^1_\ell = {\cA}_\ell
\]

\[
{\cA}^2_\ell = \{\bar a \in {\cA}_\ell:\text{ the number of }
\bar b \text{ such that } (\bar a,\bar b) \in \bold F^\ell_\beta
\text{ is } \le \lambda\}
\]

\[
{\cA}^3_\ell = {\cA}^2_\ell \cup \{\bar a:\text{ for some }
\bar b \in {\cA}^2_{3-\ell} \text{ we have } (\bar a,\bar b) \in 
\bold F^\ell_\beta\}.
\]

\mn
So $|{\cA}^3_\ell| \le \lambda$ by clause (d) of the assumption
 and let $\langle \bar a^\ell_\zeta:\zeta < \lambda\rangle$ list 
${\cA}^3_\ell$ possibly with repetitions
\mn
\begin{enumerate}
\item[$(*)$]  it is enough to take care of ${\cA}^3_\ell \cap
{\cA}_\ell$ for $\ell=1,2$.
\end{enumerate}
\mn
[Why?  By the basic properties of dependence relation.]

So we can continue.

Let $S$ be the set of limit ordinals $\delta < \lambda$ such that: for
a club of $\delta_* \in [\delta,\lambda)$ of cofinality $\aleph_0$ we
can find $\langle \bar b^\ell_\zeta:\zeta \in
[\delta,\delta_*)\rangle$ for $\ell=1,2$ such that:
\mn
\begin{enumerate}
\item[$(\alpha)$]   $\bar b^\ell_\zeta \in \{\bar a^\ell_\xi:\xi \in
[\delta,\delta_*)\}$ 
\sn
\item[$(\beta)$]   $(\bar b^1_\zeta,\bar b^2_\zeta) \in \bold F^1_\beta$
\sn
\item[$(\gamma)$]  $\langle \bar b^\ell_\zeta:\zeta \in
[\delta,\delta_*)\rangle$ is $\bold R_\ell$-independent over $\{\bar
a^\ell_\zeta:\zeta < \delta\}$
\sn
\item[$(\delta)$]  if $\zeta < \delta_*$ and $\bar a^\ell_\zeta \in
{\cA}_\ell$ then $\bar a^\ell_\zeta$ does $R_\ell$-depend on
$\{\bar a^\ell_\zeta:\zeta < \delta\} \cup \{\bar b^\ell_\zeta:\zeta
\in [\delta,\delta_*)\}$.
\end{enumerate}
\mn
If $S$ is not stationary we can easily finish (we start by playing $\omega$
moves in ${\cG}^\gamma_\mu$).  So assume $S$ is
stationary, hence for some regular $\sigma \le \lambda_1$ the set $S'
= \{\delta \in S:\text{ cf}(\delta) = \sigma\}$ is stationary.  By
playing $\sigma + \omega$ moves (recalling $\lambda \in \check
I[\lambda]$) we get a contradiction to the definition of $S$.

\noindent
2),3) Obvious. 
\end{PROOF}

\begin{remark}
In \ref{1a.23}(1), to get the exact 
$\gamma(*)$, we combine partial isomorphisms.  So
we simulate two plays and use the composition of the $f^{\bold
x^i_\beta}$'s from two plays where in each ISO use a winning strategy
{\bf st}.  
\end{remark}

\begin{claim}
\label{1a.22}
We can use a variant of Definition \ref{1a.9}(2) as follows: we can in
case 2 make a $\bold R_\ell$ dependence relation on $\kappa \times
{}^{\theta >}(M_\ell)$, but equivalently $C \times {}^{\theta
>}(M_\ell)$ for a set $C$ of cardinality $\le \kappa$, \underline{but}
\mn
\begin{enumerate}
\item[$(a)$]  it seems to help presently relevant only for $\kappa
\le 2^{|\tau(M_1)|+\aleph_0}$
\sn
\item[$(b)$]   if $\kappa \le 2^{< \theta}$ we get an equivalent game.
\end{enumerate}
\end{claim}

\begin{remark}
1) We can replace $2^{< \theta}$ by a larger cardinal
in clause (b) for ``interesting" cases of $M_1,M_2$.

\noindent
2) Anyhow we use only (b).
\end{remark}

\begin{PROOF}{\ref{1a.22}}
Clause (a) is obvious.

For clause (b), \wilog \, $\|M_\ell\| > 1$, now let $\langle
\eta_\alpha:\alpha < \kappa\rangle$ be a sequence of pairwise distinct
members of ${}^{\kappa >}2$.  now we define $F_\ell:{}^{\theta
>}(M_\ell) \rightarrow (2^{< \theta}) \times
{}^{\theta >}(M_\ell)$ as follows: for $\bar a
\in {}^{\theta >}(M_\ell)$ let

\[
\bold i(\bar a) = \text{ min}\{i:2i \ge \ell g(\bar a) \text{ or } 2i+1 <
\ell g(\bar a) \wedge a_{2i} \ne a_{2i+1}\}
\]

\[
\eta_{\bar a} = \langle T.V.(a_{2{\bold i}(\bar a)+2+2j} = 
a_{2{\bold i}(\bar a)+2+2j+1}):j\ge 0 \text{ and } 
2 \bold i(\bar a)+2+2j+1 < \ell g(\bar a)\rangle
\]

where T.V. stands for ``truth value".

\[
\alpha(\bar a) = \text{ Min}\{\alpha \le \kappa:\text{ if } \alpha <
\kappa \text{ then } \eta_\alpha = \eta_{\bar a}\}.
\]

\mn
Finally, $F_\ell(\bar a)$ 
\underline{is} $(\alpha(\bar a),\langle a_{2j}:j<i(\bar
a)\rangle)$ if $\bold i(a) < \ell g(\bar a) \wedge \alpha(\bar a) <
\kappa$, and is $(0,\bar a)$ if otherwise.

Let

\begin{equation*}
\begin{array}{clcr}
\bold R'_\ell := \{{\cA}:&{\cA} \subseteq {}^{\theta >}(M_\ell)
\text{ and } \{F(\bar a):\bar a \in {\cA}\} \in {\cR} \\
  &\text{ or for some } \bar a' \ne \bar a'' \in {\cA} \text{ we
  have } F(\bar a') = \bar a''\}.
\end{array}
\end{equation*}

\mn
Now check. 
\end{PROOF}

\begin{claim}
\label{1a.24}
$M_1,M_2$ are {\rm EF}$^+_{\gamma,\theta,\mu,\lambda}$-equivalent
\when \,:
\mn
\begin{enumerate}
\item[$(a)$]   $K$ a class of $\tau_0$-structures and $\Phi \in
\Upsilon[K]$, see \ref{2b.1.7}(8), used here for $K = K_{\text{or}}
=$ the class of linear orders and $K_{\text{oi}}$, see Definition \ref{2b.1}
\sn
\item[$(b)$]   the structures $I_1,I_2 \in K$ are
{\rm EF}$^+_{\gamma,\theta,\mu,\lambda}$-equivalent
\sn
\item[$(c)$]   $M_\ell = \text{\rm EM}_\tau(I_\ell,\Phi)$ for
$\ell=1,2$ for some $\tau \subseteq \tau_\Phi$
\sn
\item[$(d)$]   $\mu \ge \aleph_0$ and $|\tau_\Phi| < \theta$.
\end{enumerate}
\end{claim}

\begin{PROOF}{\ref{1a.24}}
Let {\bf St} be a winning strategy of the ISO player in
the game ${\cG}^+_{\gamma,\theta,\mu,\lambda}(I_1,I_2)$.  We define
a strategy {\bf st}$_*$ of the ISO player in the game 
${\cG}^+_{\gamma,\theta,\mu,\lambda}(M_1,M_2)$ as follows.

During a play of it after $\beta$ moves a partial isomorphism
$f^*_\alpha$ from $M_1$ to $M_2$ has been chosen, 
\underline{but} the ISO player also simulates a
play of ${\cG}^+_{\gamma,\theta,\mu,\lambda}(I_1,I_2)$ in which we
call the function $h_\alpha$, and in which he uses the winning
strategy {\bf st} and
\mn
\begin{enumerate}
\item[$\boxplus$]   $f_\alpha \subseteq \hat h_\alpha$ where $\hat
h_\alpha$ is defined by

$\hat h_\alpha(\sigma^{M_1}(a_{t_0},\dotsc,a_{t_{n-1}})) =
\sigma^{M_\ell}(a_{h_\alpha(t_0)},\dotsc,a_{h_\alpha(t_{n-1})})$ for
$n < \omega,\sigma(x_0,\dotsc,x_{n-1})$ a term of $\tau_\Phi$ and
$t_0,\dotsc,t_{n-1} \in \text{ Dom}(h_\alpha)$.
\end{enumerate}
\mn
Why can the player ISO 
carry this strategy {\bf st}$_*$?  Suppose we arrive to the
$\beta$-th move.  The point to check is Case 2 in Definition
\ref{1a.9}(2), so the AIS player has chosen $\bold R_1,\bold
R_2,{\cA}_1,{\cA}_2$ as there.

Let $\{\bar \sigma_\zeta(\bar x_\zeta):\zeta < 2^{< \theta}\rangle$
list $\{\bar\sigma(\bar x):\bar \sigma(\bar x) = \langle \sigma_i(\bar
x):i < \ell g(\bar \sigma)\rangle,\ell g(\bar \sigma) < \theta,\ell
g(\bar x) < \theta$ and each $\sigma_i$ is a $\tau_K$-term.

Clearly ${}^{\theta >}(M_\ell) = \{\bar \sigma^{M_\ell}_\zeta(\bar
t):\zeta < 2^{< \theta}$ and $\bar t \in {}^{\ell g(\bar
x_\zeta)}(I_\ell)\}$,  so by clause (b) of \ref{1a.22}, we can assume
``$\bold R_\ell$ is a dependence relation on $\{(\zeta,\bar
t_\zeta):\zeta < 2^{< \theta},\bar t_\zeta \in {}^\theta(I_\theta)$
and $\ell g(\bar t) = \ell g(\bar t_\zeta)\}$.

That is

\begin{equation*}
\begin{array}{clcr}
\bold R'_\ell = \{u:&\{\sigma^{M_\ell}_\zeta(\bar t):(\zeta,\bar t)
\in u\} \in \bold R_\ell\} \text{ \underline{or} there are } 
(\zeta_1,\bar t_1) \ne (\zeta_2,\bar t_2)  \\
  &\text{ from } u \text{ such that } \sigma^{M_\ell}_{\zeta_1}
(\bar a_{\bar t_1}) = \sigma^{M_\ell}_{\zeta_2}(\bar a_{\bar t_2})\}.
\end{array}
\end{equation*}

\mn
The rest should be clear.
\end{PROOF}
\newpage

\section {The properties of $T$ and relevant indiscernibility} 

In \cite[Ch.VIII]{Sh:c}, \cite[Ch.VI]{Sh:e} we use as indiscernible sets trees
with $\omega +1$ levels, suitable for dealing with unsuperstable
(complete first order) theories.

Here instead we use a linear order and family of $\omega$-sequences
from it, suitable for our case.  The letters ``oi" stands for order +
increasing ($\omega$-sequences).

\begin{definition}
\label{2b.1}
1) $K^{\text{oi}}_\lambda$ is the 
class of structures  ${\bold J}$ of the form $(J,Q,P<,F_n)_{n<\omega}
=(|\bold J|,P^{\bold J},Q^{\bold J},
<^{\bold J},F_n^{\bold J})$, where $J = |\bold J|$ is a set of cardinality 
$\lambda,<^{\bold J}$ a linear order on $Q^{\bold J} \subseteq J,
P^{\bold J}=|{\bold J}| \setminus Q^{\bold J},
F^{\bold J}_n$ a unary function,
$F^{\bold J}_n \restriction Q^{\bold J}=$ the identity and 
$a \in J \setminus Q^{\bold I} \Rightarrow F^{\bold J}_n (a) \in 
Q^{\bold J}$ and $n \ne m \Rightarrow F^{\bold J}_n(a) \ne 
F^{\bold J}_m(a)$ and for simplicity
$a \ne b \in P^M \Rightarrow \bigvee\limits_{n<\omega} F_n (a) \ne F_n
(b)$; lastly,
we add $n < m \Rightarrow F^{\bold J}_n(a) <^{\bold J} F^{\bold J}_m(a)$ 
(there is a small price).  
We stipulate $F_\omega^{\bold J}=$ the identity on $|{\bold J}|$ and
$I^{\bold J} = (Q^{\bold J},<^{\bold J})$.

\noindent
1A) $K_{\text{oi}} = \cup\{K^{\text{oi}}_\lambda:\lambda$ a cardinal$\}$.

\noindent
2) For a linear order $I$ and ${\gS} \subseteq 
\text{ inc}({}^{\omega}I)$ (see Definition \ref{0z.37}(3)),
we let ${\bold J}={\bold J}_{I,{\gS}}$ be the derived member 
of $K_{\text{oi}}$ which means: $|{\bold J}| = I \cup {\gS}, 
(Q^{|{\bold J}|},<^{\bold J})=I,F^{\bold J}_n (\eta)= \eta(n)$ 
for $n<\omega, F^{\bold J}_n (t)=t$ for $t \in I$.

\noindent
3) $K^{\text{oi}}_\lambda$ is the class of linear order of cardinality
$\lambda,K_{\text{or}} = \cup\{K^{\text{or}}_\lambda:\lambda$ a cardinal$\}$.
\end{definition}

\begin{definition}
\label{2b.1.7}
1) For a vocabulary $\tau_1$ let
$\Upsilon^{\text{oi}}_{\tau_1}$ be the class of functions $\Phi$ with
domain $\{\text{tp}_{\text{qf}}(\bar t,\emptyset,\bold J):\bar t \in
{}^{\omega >}|\bold J|,\bold J \in K^{\text{oi}}\}$, see \ref{0z.7}
and if $q(x_0,\dotsc,s_{m-1}) \in \text{ Dom}(\Phi)$ then
$\Phi(q)$ is a complete quantifier free $n$-type in $\bbL(\tau_1)$
with the natural compatibility functions.

\noindent
2) Let $\Upsilon^{\text{oi}}_\kappa = \{\Phi:\Phi \in
\Upsilon^{\text{oi}}_{\tau_1}$ for some vocabulary $\tau_1$ of cardinality
$\kappa\}$.

\noindent
3) For $\Phi \in \Upsilon^{\text{oi}}_\kappa$ let $\tau(\Phi) =
\tau_\Phi$ be the vocabulary $\tau_1$ such that $\Phi \in
\Upsilon^{\text{oi}}_{\tau_1}$.

\noindent
4) For $\Phi \in \Upsilon^{\text{oi}}_\kappa,\bold J \in
K^{\text{oi}}$ let EM$(\bold J,\Phi)$ be ``the" $\tau_\Phi$-model
$M_1$ generated by $\{a_t:t \in \bold J\}$ such that: $n < \omega,\bar
t \in {}^n \bold J \Rightarrow \text{ tp}_{\text{qf}}(\langle
a_{t_0},\dotsc,a_{t_{n-1}}\rangle,\emptyset,M_1) =
\Phi(\text{tp}_{\text{qf}}(\langle t_0,
\dotsc,t_{n-1}\rangle,\emptyset,\bold J)$. 

\noindent
5) If $\tau \subseteq \tau_\Phi$ then EM$_\tau(\bold J,\Phi)$ is the
$\tau$-reduct of EM$(J,\Phi)$.

\noindent
6) Let $\Upsilon^{\text{or}}_{\tau_1},\Upsilon^{\text{or}}_\kappa$ and
   EM$(I,\Phi)$, EM$_\tau(I,\Phi)$ be defined similarly for $\bold J$
   a linear order.

\noindent
7) Let
$\Upsilon^{\omega{\text{-tr}}}_{\tau_1},\Upsilon^{\omega{\text{-tr}}}_\kappa$
and EM$(I,\Phi)$, EM$_\tau(I,\Phi)$ be defined similarly for $\bold J
\in K^\omega_{\text{tr}}$, i.e. trees with $\omega +1$ levels (with a
linear order on the successor of any member of level $< \omega$).

\noindent
8) We can above replace $K_{\text{oi}}$ by any class $K$ of
$\tau_K$-structures. 
\end{definition}

\begin{definition}
\label{2b.2A}
1) A (complete first order) $T$ is 
$\aleph_0$-independent $\equiv$ not strongly dependent (this is
from \cite[\S3]{Sh:783}, see \cite[\S1]{Sh:863}) \when \,: 
there is a sequence $\bar \varphi=\langle \varphi_n(x,\bar y_n):n <
\omega\rangle$, (may use finite $\bar x$, as usual does not matter by
\cite[2.1]{Sh:863}) of (first order)
formulas such that $T$ is consist with $\Gamma_\lambda$ for some 
($\equiv$ every $\lambda \ge \aleph_0$) where

\[
\Gamma_\lambda =
\{\varphi_n(x_\eta,\bar y^n_\alpha)^{\text{if}(\alpha=\eta(n))}
:\eta\in {}^\omega\lambda, \alpha<\lambda,n<\omega\}.
\]

\mn
2) $T$ is strongly stable \when \, it is stable and strongly dependent.
\end{definition}

\begin{claim}
\label{2b.3}
If $T$ is first order complete, 
$T_1\supseteq T$ is first order complete, without loss of
generality with Skolem
functions and $T$ is not strongly dependent \then \, we can find 
$\bar \varphi=\langle \varphi_n (x,\bar y_n):n<\omega\rangle, 
\bar y_n \trianglelefteq \bar y_{n+1}$ 
and $\varphi_n(x,\bar y_n) \in \bbL(\tau_T)$ for $n < \omega$ such that
\mn
\begin{enumerate}
\item[$\circledast$]   for any $\bold J \in K_{\text{oi}}$ we can
find $M,\langle \bar a_t:t \in \bold J\rangle$ such that
\begin{enumerate}
\item[$(a)$]   $M$ is the Skolem hull of $\{\bar a_t:t \in \bold J\}$
\sn
\item[$(b)$]    $\bar a_t \in {}^{\omega}M$ for $t \in I^{\bold J},
\bar a_\eta = \langle a_\eta\rangle \in M_1$ for $\eta \in P^{\bold J}$
\sn
\item[$(c)$]    for $\eta \in P^{\bold J},t \in Q^{\bold J}$ and
$n < \omega$ we have $M \models \varphi_n [\bar a_\eta,\bar a_t]$ 
iff $F_n(\eta)=t$; (pedantically we should write 
$\varphi_n (a_\eta,\bar a_t \restriction \ell g(\bar y_n)))$
\sn
\item[$(d)$]    $\langle \bar a_t:t \in \bold J\rangle$ is
indiscernible in $M$ for the index model $\bold J$
\sn
\item[$(e)$]    $M$ is a model of $T_1$
\sn
\item[$(f)$]   in fact (not actually used, see \ref{2b.1.7}) 
there is $\Phi \in \Upsilon^{\text{oi}}_{|T_1|}$ 
depending on $T_1,\bar\varphi$ only
such that $M = \text{\rm EM}(\bold J,\Phi)$, in fact if $n < \omega,\bar t =
\langle t_\ell:\ell < n\rangle \in \bold J$ then {\rm tp}$_{\text{qf}}(\bar
a_{t_0} \char 94 \ldots \char 94 \bar a_{t_{n-1}},\emptyset,M) =
\Phi(\text{\rm tp}_{\text{qf}}(\bar t,\emptyset,\bold J))$.
\end{enumerate}
\end{enumerate}
\end{claim}

\begin{PROOF}{\ref{2b.3}}
Let $I = (Q^{\bold J},<_{\bold J})$.
By an assumption, i.e. \ref{2b.2A} there is a sequence $\langle
\varphi'_n(x,\bar y_n):n < \omega\rangle$ as in Definition
\ref{2b.2A} and let $k_n = \ell g(\bar y_n)$.

Let $I$ be an infinite linear order.   Easily we 
can find $M_1 \models T_1$ and a sequence 
$\langle \bar a_t:t \in I \rangle$ with $\bar a_t \in 
{}^\omega(M_1)$ such that for every 
$\eta\in {}^{\omega}I$, the set 
$\{\varphi_n (\bar x,\bar a_t)^{\text{if}(\eta(n)=t)}:
t \in I,n<\omega\}$ is a type, i.e. finitely satisfiable in $M_1$. 

Now by Ramsey theorem \wilog \, $\langle \bar a_t:t \in I\rangle$ is an
indiscernible sequence in $M_1$.  Without loss of generality 
$M_1$ is $\lambda^+$-saturated, we then expand $M_1$ to $M_1^+$ by function 
$F_n^{M^+_1} (n<\omega)$, (of finite arity) such that for $t_0
<_{\bold J} \ldots <_{\bold J} t_{n-1}$ from $Q^{\bold J}$ the element
$F_n (\bar a_{t_0},\bar a_{t_1},\ldots\bar a_{t_{n-1}})$ or more
exactly $F_n(\bar a_{t_0}\restriction k_0,\bar a_{t_1} \restriction
k_1,\ldots,\bar a_{t_{n-1}} \restriction k_{n-1})$ realizes in 
$M_1$ the type $\{\varphi_\ell (x,\bar a_t)^{\text{if}(\eta(\ell)=t)}:
t \in I,\ell<n\}$.
Let $M^+_2$ be an expansion of $M^+_1$ by Skolem functions such that
$|\tau_{M^+_2}| = |T_1|$, (natural, though not strictly required).
Without loss of generality $\langle \bar a_t:t \in I \rangle$ is an
indiscernible sequence also in $M^+_2$.

Let $D$ be a non-principal ultrafilter
on $\omega$ and in $M^+_3 = (M^+_2)^\omega /D$, we let 
$\bar a'_t=\langle \bar a_t:n<\omega\rangle/D$ for $t \in I$, and 
$\bar a'_\eta=\langle F_n(\bar a_{\eta(0)},\bar a_{\eta(1)},
\ldots,\bar a_{\eta(n-1)}):n<\omega\rangle / D$ for 
$\eta \in \text{ incr}({}^{\omega}I)$ and $\bar a'_t = \bar
a'_{<F^{\bold J}_n(t):n <\omega>}$ for $t \in P^{\bold J}$.

Let $M^+_4$ be the submodel of $M^+_3$ generated by $\{\bar a'_t:t \in
\bold J\}$ and $M$ be $M^+_4 \rest \tau(T_1)$.  Now $M,\langle \bar
a_t:t \in \bold J\rangle$ are as required.
\end{PROOF}

\begin{claim}
\label{2b.4}
Assume ${\bold J}_\ell \in K_{\text{oi}}$, and 
$M_\ell,\bar \varphi,T_1,T$ as in \ref{2b.3} for
$\ell=1,2$.   A sufficient condition for 
$M_1 \restriction \tau_T \ncong M_2 \restriction \tau_T$ is:
\mn
\begin{enumerate}
\item[$(*)$]   if $f$ is a function from ${\bold J}_1$ (i.e. its
  universe) into ${\cM}_{|T_1|,\aleph_0} ({\bold J}_2)$ 
(i.e. the free algebra generated by $\{x_t:t\in {\bold J}_2\}$ in the 
vocabulary $\tau_{|T_1|,\aleph_0}= \{F^n_\alpha:n<\omega$ and 
$\alpha<|T_1|\},F^n_\alpha$ has arity $n$, see more in
\cite[Ch.III,\S1]{Sh:e} = \cite{Sh:E59})
we can find $t\in P^{{\bold J}_1}, n<\omega$, and $s_1,s_2 \in 
Q^{{\bold J}_1}$ and $k,\sigma,r^\ell_i(\ell=1,2$ and $i<k),m,\sigma^*$
such that:
\begin{enumerate}
\item[$(\alpha)$]   $F^{{\bold J}_1}_n (t)= s_1 \ne s_2$
\sn
\item[$(\beta)$]   for $\ell \in \{1,2\}$ we have
$f(s_\ell)= \sigma (r^\ell_0,\ldots, r^\ell_{k-1})$ 
so $k<\omega,r_t^\ell\in {\bold J}_2$ for $i<k$ and 
$\sigma$ is a $\tau_{|T_1|,\aleph_0}$-term
not dependent on $\ell$
\sn
\item[$(\gamma)$]   $f(t)=\sigma^* (r_0,\ldots, r_{m-1}),
\sigma^*$ is a $\tau_{|T_1|,\aleph_0}$-term and  
$r_0,\ldots,r_{m-1}\in {\bold J}_2$
\sn
\item[$(\delta)$]    the sequences 
\[
\langle r^1_i:i<k\rangle \char 94 \langle r_i:i<m\rangle
\]

\[
\langle r^2_i:i<k\rangle \char 94 \langle r_i:i<m\rangle
\]

realize the same quantifier free type in ${\bold J}_2$
(note: we should close by the $F^{{\bold J}_2}_n$'s, but no need to
iterate as $F^{\bold J_2}_n \rest Q^{\bold J_2}$ is the identity so
quantifier free 
type mean the truth value of the inequalities $F_{n_1} (r')\ne
F_{n_2}(r')$ (including $F_\omega$) and the order between those terms).
\end{enumerate}
\end{enumerate}
\end{claim}

\begin{PROOF}{\ref{2b.4}}
Straight (or as in \cite[Ch.III]{Sh:e} = \cite{Sh:E59}).
\end{PROOF}

\begin{remark}
We could have replaced $Q^{\bold J}$ by the disjoint 
union of $\langle Q^{\bold J}_n:n<\omega\rangle,<^{\bold J}$ linearly 
order each $Q^{{\bold J}}_n$ (and $<^{\bold J} = \cup 
\{< \restriction Q^{{\bold J}_1}_n:n<\omega\}$ and use 
$Q_n$ to index parameters for $\varphi_n (x,\bar y_n))$.
Does not matter at present. 
\end{remark}

But for our aim we can replace ``not strongly stable" by a weaker
demand, though this will not be carried here, we present it.
Recall (from \cite[\S5(G)]{Sh:863},
i.e. \cite[5.39=dw5.35tex(2A)]{Sh:863}) the following, an equivalent
definition to ``a (complete first order theory) $T$ is 
strongly$_4$ dependent".  
\begin{definition}
\label{2b.5}
1)  A (complete first order) $T$ is
 not strongly$_4$ dependent \If \, there is a sequence 
$\bar \varphi=\langle \varphi_n (\bar x,\bar y_n):n<\omega\rangle$, 
(finite $\bar x$ of length $m < \omega$, as usual) of 
(first order) formulas from $\bbL(\tau_T)$, 
an infinite linear order $I$, a sequence
$\langle \bar a_\eta:\eta \in \text{ incr}_{< \omega}(I)\rangle$
indiscernible in $M$ with $\ell
g(\bar a_\eta) \le \omega$ and letting $B = \cup\{\bar a_\eta:\eta \in
\text{ incr}_{< \omega}(I)\rangle$ for some $m < \omega$ and 
$p \in \bold S^m(B,M)$ for every $k < \omega$ there is $n < \omega$,
\underline{satisfying}: for no linear order $I^+$ extending $I$ and subset
$I_0$ of $I^+$ with $\le k$ members, do we have:
\mn
\begin{enumerate}
\item[$\otimes$]   if $\bar t^1,\bar t^2$ are increasing sequences from
$I$ of the same length $n$ realizing the same quantifier-free
type over $I_0$ in $I^+$ and for $i=1,2$ we let
$\bar b^i = (\ldots \char 94(\bar a_{\langle t^i_{\eta(\ell)}:
\ell < \ell g(\eta)\rangle} \restriction n) \char 94 \ldots)_{\eta \in 
\text{ inc}_{< \omega}(n)}$ \then \,
$\ell < n  \wedge u \subseteq \ell g(\bar b^1)
\wedge |u| = \ell g(\bar y_\ell) \Rightarrow \varphi_\ell(x,\bar b^1
\restriction u) \in p \Leftrightarrow \varphi_\ell(\bar x,\bar b^2
\restriction u) \in p$.
\end{enumerate}
\mn
1A) In (1) \wilog \, $\bar y_n \triangleleft \bar y_{n+1}$ for $n <
\omega$.

\noindent
2) $T$ is strongly$_4$ stable if it is stable and strongly dependent.
\end{definition}

\begin{remark}
1) We can write the condition in \ref{2b.5}(1)
 without $I^+$ speaking on finite sets as done in $(*)$
in the proof of \ref{2b.6} below.

\noindent
2) In \ref{2b.5} by compactness we can get such
$\langle \bar a'_\rho:\rho \in \text{ inc}_{< \omega}(I')\rangle$ for
any infinite linear order $I'$.  

Next we deduce a consequence of being non-strongly$_4$-dependent, see 
Definition \ref{2b.5} helpful in proving non-structure results.
\end{remark}

\begin{claim}
\label{2b.6}
If $T$ is first order complete, 
$T_1\supseteq T$ is first order complete, without loss of generality 
with Skolem functions and 
$T$ is not strongly$_4$ dependent as witnessed by $\bar\varphi =
\langle \varphi_n(\bar x,\bar y_n):n<\omega\rangle$, i.e. as in Definition
\ref{2b.5}(1A), \then \, there is $\tau_1 \supseteq
\tau_{T_1},|\tau_1| = |T_1|$ and
$\bar \sigma_n(\bar z_n) = \langle \sigma_{n,\ell}(\bar z_n):
\ell < \ell g(\bar\sigma_n)\rangle,\sigma_{n,\ell}$ is a $\tau_1$-term
such that:
\mn
\begin{enumerate}
\item[$\circledast$]    if $I,{\gS}$ and $\bold J = \bold
J_{I,{\gS}}$ are as in Definition \ref{2b.1}(2),
\then \, there are $M_1$ and $\langle \bar a_t:t\in I\rangle$ 
and $\langle \bar a_\eta:\eta\in {\gS}\rangle$ such that:
\begin{enumerate}
\item[$(\alpha)$]    $M_1$ is a $\tau_1$-model and is
the Skolem hull of $\{\bar a_t:t\in I\} \cup 
\{\bar a_\eta:\eta \in {\gS}\}$

(we write $\bar a_t$ for $t \in {\gS} \subseteq \bold J$ for uniformity)
\sn
\item[$(\beta)$]  $\langle \bar a_t:t \in \bold J\rangle$ is 
indiscernible in $M_1$,
\sn
\item[$(\gamma)$]     if $\eta \in {\gS}$ and $k < \omega$
\then \, for large enough $n(*)$ we have:
\sn
\item[${{}}$]   $(*) \quad$ if $u \subseteq n(*),|u| \le k$, then 
we can find $\bar s,\bar t$ and $n_* < n(*)$ and $\bar \sigma$
such that
\sn
\item[${{}}$]  $\quad (i) \quad \bar s,\bar t$ are sequences of members of
$\{F^{\bold J}_n(\eta):n < n(*)\}$
\sn
\item[${{}}$]  $\quad (ii) \quad \ell g(\bar s) = \ell g(\bar t) \le n(*)$
\sn
\item[${{}}$]  $\quad (iii) \quad s_i <_I s_j 
\Leftrightarrow t_i <_I t_j$ for $i,j < \ell g(\bar s)$
\sn
\item[${{}}$]  $\quad (iv) \quad$ if $i < \ell g(\bar s) = \ell g(\bar t)$ then

\hskip35pt $(\forall n \in u)(F^{\bold J}_n(\eta) \le_I s_i \equiv 
F^{\bold J}_n(\eta) \le_I t_i)$
\sn
\item[${{}}$]  $\quad (v) \quad \bar \sigma = 
\langle \sigma_i(\bar y):i < \ell
g(\bar y_{n_*})\rangle,\sigma_i$ a $\tau_1$-term
\sn
\item[${{}}$]  $\quad (vi) \quad M_1 \models \varphi_{n_*}[\bar a_\eta,\dotsc,
\sigma^{M_1}_i(\bar a_{t_0},\bar a_{t_1},\ldots),\ldots]_{i < \ell g
(\bar y_{n_*})} \equiv$

\hskip35pt $\neg \varphi_{n_*}[\bar a_\eta,\dotsc,
\sigma^{M_1}_i(\bar a_{s_0},\bar a_{s_1},\ldots),
\ldots]_{i < \ell g(\bar y_{n_*})}$ 
\sn
\item[$(\delta)$]    $M_1$ is a model of $T_1$ so $\tau_{M_1}
\supseteq \tau_{T_1}$.
\end{enumerate}
\end{enumerate}
\end{claim}

\begin{PROOF}{\ref{2b.6}}  
Fix $I,{\gS}$; \wilog \, $I$ is dense with neither
first nor last element and is $\aleph_1$-homogeneous hence
there are infinite increasing sequences of members of $I$.

Let $I,\langle \varphi_n(\bar x,\bar y_n):n < \omega\rangle,\langle \bar
a_\eta:\eta \in \text{ incr}_{< \omega}(I)\rangle$
 and $p \in \bold S^m(\cup\{\bar a_\eta:\eta \in \text{
incr}_{<\omega}(I)\})$ exemplify $T$ is not strongly$_4$ dependent, i.e.
be as in the Definition so $m = \ell g(\bar x)$.  For notational
simplicity (and even \wilog \, by \cite[\S5]{Sh:863}) assume $m=1$.

Now in \ref{2b.5} we can add:
\mn
\begin{enumerate}
\item[$(*)$]  there is a sequence $\langle(n_k,m_k,I^*_k):k <
\omega\rangle$ such that $k < n_k < m_k,m_k < m_{k+1},I^*_k \subseteq
I$ has $m_k$ members, for no $I_0 \subseteq I^*_k$ with $\le k$
members does $\otimes$ from \ref{2b.5} holds for $\bar t^1,\bar t^2
\in \text{ incr}_{< n_k}(I^*_k)$ and $(k,n_k)$ here standing for $k,n$.
\end{enumerate}
\mn
[Why?  By compactness.]

Without loss of generality $I$ is the reduct to the vocabulary
$\{<\}$, i.e. to just a linear order of an ordered
field $\bbF$ and $t_q \in \bbF$ for $q \in \bbQ$ are such
that $0 <_{\bbF} t_q,(t_{q_1})^2 <_{\bbF} t_{q_2}$ for $q_1 <_{\bbF} q_2$
(hence $n < \omega \Rightarrow n <_{\bbF} t^n_{q_1} <_{\bbF}
t_{q_2}$).  By easy manipulation \wilog \, $I^*_k 
= \{t_i:i=0,1,\dotsc,m_k\}$.  

Now for each $m < \omega$ and $\eta \in \text{ incr}_m(I)$ we can
choose $c_\eta$ such that if $m = m_k$ then for some automorphism $h$
of $I$ mapping $I^*_k$ onto Rang$(\eta)$, letting $\hat h$ be an
automorphism of $M_1$ mapping $\bar a_\nu$ to $\bar a_{h(\nu)}$
for $\nu \in \text{ incr}_{< \omega}(I)$, the element
$c_\eta$ realizes $\hat h(p)$ and $\langle c_\eta:\eta \in \text{ incr}_{<
\omega}(I)\rangle$ is without repetitions.

Now \wilog \, $\left< \langle c_\eta\rangle \char 94 \bar a_\eta:\eta
\in \text{ incr}_{< \omega}(I)\right>$ is an indiscernible sequence
and let $a_t = c_{<t>}$ be such that $M_0$ be a model of $T_1$ satisfying
$\cup\{\langle c_\eta\rangle \char 94 \bar a_\eta:\eta \in \text{
incr}_{< \omega}(I)\} \subseteq M_0 \restriction \tau \prec {\gC}$.  
Without loss of generality $\left< \langle c_\eta \rangle \char
94 \bar a_\eta:\eta \in \text{ incr}_{< \omega}(I)\right>$ is
indiscernible in $M_0$ and we can find an expansion $M_1$ of $M_0$
such that $|\tau_{M_2}| = |T_1|$ such that $\bar a_\eta = \langle
F_{\ell g(\eta),i}(\bar a_{\eta(0)},\dotsc,\bar a_{\eta(n-1)}):i < \ell g(\bar
a_\eta)\rangle,c_\eta = F_{\ell g(\eta)}(\bar a_{\eta(0)},
\dotsc,\bar a_{\eta(n-1)})$ if $\eta \in 
\text{ incr}_n(I)$ and $M_1$ has Skolem functions.

By manipulating $I$ \wilog \, we can find $I_* \subseteq I$ of order
type $\omega$. 

So for some $H_n \in \tau_1$ for $n < \omega$,
\mn
\begin{enumerate}
\item[$\odot$]   if $t_0 < t_1 < \ldots$ list $I_*$, for every $k <
\omega$ large enough, for every $u \subseteq n(*)$ satisfying
$|u| \le k$ 
for every $n$ large enough $H^{M_1}_n(\bar a_{t_0},
\bar a_{t_1},\dotsc,\bar a_{t_{n-1}})$ satisfies
the demand (on the singleton $\bar a_\eta$ 
from clause $(\gamma)$ in the claim).
\end{enumerate}
\mn
Let $D$ be a non-principal ultrafilter on $\omega$ such that
$\{m_k:k < \omega\} \in D$, let $M_2$ be
isomorphic to $M^\omega_1/D$ over $M_1$, i.e. $M_1 \prec M_2$
and there is an isomorphism $\bold f$ from $M_2$ onto
$M^\omega_1/D$ extending the canonical embedding.

If $\eta$ is an increasing $\omega$-sequence of members of $I$, we let

\[
a^n_\eta = H^{M_1}_n(a_{\eta(0)},\dotsc,a_{\eta(n-1)}) \in M_1
\]

\mn
and let

\[
a_\eta = \bold j^{-1}(\langle
a^0_\eta,a^1_\eta,\dotsc,a^n_\eta,\dotsc:n < \omega\rangle/D) \in M_2.
\]

\mn
Let $M'_2$ be the Skolem hull of $\{\bar a_t:t \in I\} \cup \{a_\eta:\eta
\in {\gS}\}$ inside $M_2$.  It is easy to check that it is as required.
\end{PROOF}

\noindent
Naturally it is helpful to have a sufficient condition for the
non-isomorphism of two such models:
\begin{claim}
\label{2b.7}
Assume ${\bold J}_\ell \in K^{\text{oi}}$, and 
$M_\ell,\bar \varphi,T_1,T$ as in \ref{2b.6} for $\ell=1,2$. 
A sufficient condition for $M_1 \ncong M_2$ is 
\mn
\begin{enumerate}
\item[$(*)$]    if $f$ is a function from ${\bold J}_1$ 
(i.e. its universe) into ${\cM}_{|T_1|,\aleph_0} ({\bold J}_2)$ 
(i.e. the free algebra generated by $\{x_t:t\in {\bold J}_1\}$ the 
vocabulary $\tau_{|T_1|,\aleph_0}= \{F^n_\alpha:n<\omega$ and 
$\alpha<|T_1|\},F^n_\alpha$ has arity $n$),
we can find $t\in P^{{\bold J}_1}$ and $k_* < \omega$ such that for
every $n_* < \omega$ we can find $\bar s_1,\bar s_2$ such that:
\begin{enumerate}
\item[$(\alpha)$]   $\bar s_1,\bar s_2 \in^k {}^{n_*}(Q^{\bold J})$ are
increasing, $\bar s_1 = \langle F^{\bold J}_n(t):n < n_*\rangle$ and $n <
k_* \Rightarrow s_{2,n} = s_{1,n}$ and $s_{1,n_*-1} <_I s_{2,k_*}$
\sn
\item[$(\beta)$]   $f(\bar s_\ell)= \bar \sigma(r^\ell_0,\ldots,
r^\ell_{k-1})$ so $k<\omega,r_t^\ell\in {\bold J}_2$ for 
$i<k$ so $\sigma$ is a $\tau_{|T_1|,\aleph_0}$-term
not dependent on $\ell$
\sn
\item[$(\gamma)$]    $f(t)=\sigma^* (r_0,\ldots,r_{m-1}),
\sigma^*$ is a $\tau_{|T_1|,\aleph_0}$-term and 
$r_0,\ldots,r_{m-1}\in {\bold J}_2$
\sn
\item[$(\delta)$]  the sequences 
\[
\langle r^1_i:i<k\rangle \char 94 \langle r_i:i<m\rangle
\]

\[
\langle r^2_i:i<k\rangle \char 94 \langle r_i:i<m\rangle
\]
realize the same quantifier free type in ${\bold J}_2$
(note: we should close by the $F^{{\bold J}_2}_n$, so 
type mean the truth value of the inequalities $F_{n_1} (r') \ne
F_{n_2}(r')$ (including $F_\omega$) and the order between those 
terms).
\end{enumerate}
\end{enumerate}
\end{claim}

\begin{PROOF}{\ref{2b.7}}
As in \cite[Ch.III]{Sh:300} or better in
\cite[Ch.III]{Sh:e} = \cite{Sh:E59}, called unembeddability.  
\end{PROOF}
\newpage

\section {Forcing EF$^+$-equivalent Consistency non-isomorphic models} 

The following result is not optimal, but it is enough to prove
necessary conditions on $T$ for being lean and even on $(T,*)$.
As for unstable $T$, see below in \S4.
\sn
So our main result is
\begin{claim}
\label{3c.4}
Assume $(\bar \varphi,T,T_1,\Phi)$ is 
as in \ref{2b.3}, $T$ stable and $\lambda=\lambda^{<\lambda} \ge
\aleph_1 + |T_1|$ and $\mu = \lambda^+ > \lambda$.  
\Then \, for some $\lambda$-complete 
$\lambda^+$-c.c. forcing notion ${\Bbb Q}$ we have:
$\Vdash_{\bbQ}$ ``there are models $M_1,M_2$ of $T$
of cardinality $\lambda^+$ such that $M_1 \restriction \tau(T),M_2
\restriction \tau(T)$ are {\rm EF}$^+_{\alpha,\lambda,\lambda^+}$-equivalent
for every $\alpha<\lambda$ but are not isomorphic".
\end{claim}

\begin{remark}
\label{3c.49}
1) It should be clear that we can improve it 
allowing $\alpha<\lambda^+$ and replacing forcing by e.g. 
$2^\lambda = \lambda^+$ and $\lambda= \lambda^{<\lambda}$, but we shall
continue in \cite{Sh:F918}.
\end{remark}

\begin{PROOF}{\ref{3c.4}}
We define ${\bbQ}$ as follows:
\mn
\begin{enumerate}
\item[$\circledast_1$]    $p \in {\bbQ}$ iff $p$ 
consist of the following objects satisfying the following conditions:
\begin{enumerate}
\item[$(a)$]   $u=u^p \in [\mu]^{<\lambda}$ such that 
$\alpha + i \in u\wedge i<\lambda \Rightarrow \alpha\in u$
\sn
\item[$(b)$] $<^p$ a linear order of $u$ such that
\[
\alpha,\beta\in u \wedge 
\alpha+\lambda \le \beta \Rightarrow \alpha<^p \beta 
\]
\sn
\item[$(c)$]  for $\ell=1,2\quad {\frak S}^p_\ell$ 
is a subset of $\{\eta\in {}^{\omega}u:\eta(n)+\lambda\le \eta 
(n+1)$ for $n<\omega\}$ such that $\eta \ne \nu \in 
{\gS}^p_\ell \Rightarrow \text{\rm Rang}(\eta) \cap \text{\rm
  Rang}(\nu)$ is finite; note that in 
particular $\eta \in {\gS}^p_\ell$ is without repetitions and is
$<^p$-increasing 
\sn
\item[$(d)$]  $\Lambda^p$ a set of $<\lambda$ increasing sequences 
of ordinals 
from $\{\alpha\in u^p:\lambda| \alpha\}$ hence of length $<\lambda$
\sn
\item[$(e)$]   $\bar f^p = \langle f^p_\rho:\rho\in\Lambda^p\rangle$ 

such that
\sn
\item[$(f)$]   $f^p_\rho$ is a partial automorphism of 
the linear order $(u^p,<^p)$ such that $\alpha \in 
\text{ Dom}(f^p_\rho) \Rightarrow \alpha 
+ \lambda = f^p_\rho(\alpha) + \lambda$
 and we let $f^{1,p}_\rho = f^p_\rho,f^{2,p}_\rho=(f^p_\rho)^{-1}$
\sn
\item[$(g)$]  if $\eta \in {\gS}^p_\ell, \rho\in 
\Lambda^p, \ell\in \{1,2\}$ then {\rm Rang}$(\eta)$ is 
included in {\rm Dom}$(f^{\ell,p}_\rho)$ or is almost 
disjoint to it (i.e. except finitely many ``errors'')
\sn
\item[$(h)$]   if $\rho \triangleleft \varrho \in \Lambda^p$ 
then $\rho \in \Lambda^p$ and $f^p_\rho \subseteq f^p_\varrho$ 
\sn
\item[$(i)$]  $f^{\ell,p}_{<>}$ is the empty function and
 if $\rho \in \Lambda^p$ has limit length \then \,
\[
f^p_\rho=\cup \{f^p_{\rho\restriction i}:i < \ell g(\rho)\}
\]
\sn
\item[$(j)$]   if $\rho \in \Lambda^p$ has length $i+1$ 
then {\rm Dom}$(f^{\ell,p}_\rho)\subseteq \rho (i)$ for $\ell=1,2$
\sn
\item[$(k)$]  if $\rho \in \Lambda^p$ and $\eta\in 
{}^\omega (\text{\rm Dom}(f^p_\rho))$ then 
$\eta\in {\gS}^p_1 \Leftrightarrow 
\langle f^p_\rho (\eta(n)):n<\omega\rangle \in {\gS}^p_2$
\sn
\item[$(\ell)$] if $\rho_n \in \Lambda^p$ for $n<\omega$
and $\rho_n \triangleleft \rho_{n+1}$ and $\lambda>\aleph_0$
then $\cup\{\rho_n:n<\omega\}\in \Lambda$
\end{enumerate}
\item[$\circledast_2$]   We define the order $\le = \le_{\bbQ}$ 
on ${\bbQ}$ as follows:
$p \le q$ iff $(p,q\in {\bbQ}$ and)
\begin{enumerate}
\item[$(a)$]   $u^p\subseteq u^q$
\sn
\item[$(b)$]  $\le^p = \le^q \restriction u^p$
\sn
\item[$(c)$]   ${\gS}^p_\ell\subseteq {\gS}^q_\ell$ for $\ell=1,2$
\sn
\item[$(d)$]    $\Lambda^p\subseteq \Lambda^q$
\sn
\item[$(e)$]  if $\rho \in \Lambda^p$ then 
$f^p_\rho \subseteq f^q_\rho$ 
\sn
\item[$(f)$]  if $\eta \in {\gS}^q_\ell\setminus 
{\gS}^p_\ell$ then {\rm Rang}$(\eta) \cap u^p$ is finite 
\sn
\item[$(g)$]  if $\rho\in \Lambda^p$ and 
$f^p_\rho \neq f^q_\rho$ \then \, 
$u^p \cap \text{\rm sup Rang}(\rho)\subseteq \text{\rm Dom}
(f^{\ell,q}_\rho)$ for $\ell=1,2$
\sn
\item[$(h)$]   if $\rho\in \Lambda^p$ and $\ell\in
\{1,2\},\alpha \in u^p \setminus \text{\rm Dom}(f^{\ell,p}_\rho)$ and 
$\alpha \in \text{\rm Dom}(f^{\ell,q}_\rho)$ then 
$f^{\ell,p}_\rho (\alpha) \notin u^p$
\end{enumerate}
\end{enumerate}
\mn
Having defined the forcing notion ${\bbQ}$ we start to 
investigate it.
\mn
\begin{enumerate}
\item[$\circledast_3$]    ${\bbQ}$ is a partial order of 
cardinality $\mu^{< \lambda} = \lambda^+$.
\end{enumerate}
\mn
[Why?  Obviously.]
\mn
\begin{enumerate}
\item[$\circledast_4$]   $(i) \quad$ 
if $\bar p=\langle p_i:i<\delta\rangle$ 
is $\le^{\bbQ}$-increasing, $\delta$ a limit ordinal $<\lambda$ 
of uncountable cofinality then 
$p_\delta := \cup\{p_i:i<\delta\}$ defined naturally is an 
upper bound of $\bar p$ 
\sn
\item[${{}}$]    $(ii) \quad$ if $\delta<\lambda$ is a limit 
ordinal of cofinality $\aleph_0$ and the sequence 

\hskip25pt $\bar p=\langle p_i:i<\delta\rangle$ 
is increasing (in ${\bbQ}$), \then \, it has an upper bound. 
\end{enumerate}
\mn
Why $(i)$?  Think or see (ii); why the 
case cf$(\delta) > \aleph_0$ is easier?  
Because of clause $\circledast_1(i)$ and $\circledast_1(\ell)$.
Why $(ii)$?  We define $p_\delta \in {\bbQ}$ as follows:
$u^{p_\delta} = \cup \{u^{p_i}:i<\delta\}$,
$<^{p_\delta} = \cup \{<^{p_i}:i<\delta\},\Lambda^{p_\delta}
=\cup \{\Lambda^{p_i}:i<\delta\}\cup \{\rho:\rho$ is an increasing 
sequence of ordinals from $u^{p_\delta}$ of length 
a limit ordinal of cofinality $\aleph_0$ such that 
$\varepsilon < \ell g(\rho)\Rightarrow 
\rho \restriction \varepsilon \in \cup \{\Lambda^{p_i}:i<\delta\}\}$. 

Let $\bar f^{p_\delta} = \langle f^{p_\delta}_\rho:\rho \in
\Lambda^q\rangle$ where: if $i < \delta$ and $\rho \in \Lambda^{p_i}
\backslash \cup\{\Lambda^{p_j}:j<i\}$, then $f^q_\rho =
\cup\{f^{p_j}_\rho:j \in [i,\delta)\}$ and if $\rho \in
\Lambda^{p_\delta} \backslash \{\Lambda^{p_i}:i<\delta\}$ 
then $f^{p_\delta}_\rho
= \cup\{f^{p_\delta}_{\rho \rest \varepsilon}:\varepsilon < \ell
g(\rho)\}$ is well defined as $\varepsilon < \ell g(\rho) \Rightarrow
\rho \char 94 \langle \varepsilon \rangle 
\in \cup\{\Lambda^{p_j};j < \delta\}$.  
Clearly clauses (a),(b),(d),(e),(f),(h),(i),(j),$(\ell)$
from $\circledast_1$ for $p_\delta \in \bbQ$ hold.

Lastly, let ${\gS}^{p_\delta}_\ell = 
\cup\{{\gS}^{p_\alpha}_\ell:\alpha < \delta\}$ for $\ell=1,2$.

Note
\mn
\begin{enumerate}
\item[$\odot_1$]   if $\rho \in \Lambda^{p_\delta} \backslash
\cup\{\Lambda^{p_\alpha}:\alpha < \delta\}$ then
Dom$(f^{p_\delta}_\rho) = u^{p_\delta} \cap \text{ sup Rang}(\rho) = 
\text{ Rang}(f^{p_\delta}_\rho)$ and for every $\alpha < \delta$ for
some $\beta < \delta$ we have $f^{p_\delta}_\rho \rest u^{p_\alpha}
\subseteq f^{p_\beta}_{\rho \rest i}$ for some $i < \ell g(\rho)$.
\end{enumerate}
\mn
[Why?  Clearly, cf$(\ell g(\rho)) = \aleph_0$.

Assume $\alpha < \delta$ and $i < \ell g(\rho)$.  Clearly for some
$\beta \in (\alpha,\delta)$ we have $\rho \rest i \in
\Lambda^{p_\beta}$.  Also the set $\{j < \ell g(\rho):\rho
\rest j \in \Lambda^{p_\beta}\}$ is an initial segment of $\ell
g(\rho)$ and cannot be $\ell g(\rho)$ ecause $\rho \notin \Lambda^{p_\beta}$
by clause $\circledast_1(\ell)$.  So for some $j < \ell g(\rho)$ we have
$\rho \rest j \notin \Lambda^{p_\beta}$ but by the choice of
$\rho$ for some $\gamma < \delta$ we have $\rho \rest j \in
\Lambda^{p_\gamma}$, so necessarily $\beta < \gamma$.  As $p_\alpha
\le_{\bbQ} p_\beta \le_{\bbQ} p_\gamma$ by clause (g) of
$\circledast_2$, as $\rho \rest i \in \Lambda^{p_\gamma} 
\backslash \Lambda^{p_\beta}$ we know that $u^{p_\beta} 
\cap \text{ sup Rang}(\rho \rest i)$ is 
included in Dom$(f^{\ell,p_\gamma}_{\rho \rest i})$ for
$\ell=1,2$ by $p_\alpha \le_{\Bbb Q} p_\beta$ hence $u^{p_\alpha} \cap
\text{ sup Rang}(\rho \rest i)$ is included in
Dom$(f^{\ell,p_\gamma}_{\rho \rest i_\alpha})$ which $\subseteq$
Dom$(f^{\ell,p_\delta}_{\rho \rest i})$ for $\ell=1,2$.  

As this holds for any $\alpha < \delta$ and $i < \ell g(\rho)$
and $u^{p_\delta} \cap \text{ sup Rang}(\rho \rest i) = 
\cup\{u^{p_\alpha} \cap \text{ sup
Rang}(\rho):\alpha < \delta\}$ it follows that
for $\ell=1,2$ we have $\varepsilon \in u^{p_\delta} \cap \text{ sup
Rang}(\rho) \Rightarrow (\exists \alpha < \delta)
(\varepsilon \in u^{p_\alpha} \cap \text{ sup Rang}(\rho)) 
\Rightarrow (\exists \beta < \delta)
[\varepsilon \in \text{ Dom}(f^{\ell,p_\delta}_\rho)]
\Rightarrow \varepsilon \in \text{ Dom}(f^{\ell,p_\delta}_\rho)$ so are done.]
\mn
\begin{enumerate}
\item[$\odot_2$]   if $\rho \in \cup\{\Lambda^{p_\alpha}:\alpha <
\delta\}$ then exactly one of the following occurs:
\begin{enumerate}
\item[$(a)$]   there is a unique $\alpha = \alpha(\rho) < \delta$ such
that $\rho \in \Lambda^{p_\alpha},(\forall \beta)(\alpha \le \beta <
\delta \Rightarrow f^{p_\beta}_\rho = f^{p_\alpha}_\rho)$ and
$(\forall \beta < \alpha)(\rho \in \Lambda^{p_\beta} \rightarrow
f^{p_\beta}_\rho \ne f^{p_\alpha}_\rho)$
\sn
\item[$(b)$]  Dom$(f^{p_\delta}_\rho) = u^{p_\delta} \cap \text{
sup Rang}(\rho) = \text{ Rang}(f^{p_\delta}_\rho)$ and 
$(\forall \alpha < \delta)(\exists \beta < \delta)(f^{p_\delta}_\rho \rest
u^{p_\alpha} \subseteq f^{p_\beta}_\rho)$.
\end{enumerate}
\end{enumerate}
\mn
[Why?  Similarly to the proof of $\odot_1$.]  

To finish proving $p_\delta \in \bbQ$, i.e. verifying
$\circledast_1$ holds, we have to check clauses (c),(g),(k).
\medskip

\noindent
\underline{Clause $(c)$}:  Obvious by the choice of ${\gS}^{p_\delta}_1$.
\medskip

\noindent
\underline{Clause $(g)$}:

So let $\eta \in {\gS}^{p_\delta}_\ell,\rho \in
\Lambda^{p_\delta}$ where $\ell \in \{1,2\}$ and we should prove that
Rang$(\eta) \subseteq \text{ Dom}(f^{\ell,p,\delta}_\rho)$ or Rang$(\eta) \cap
\text{ Dom}(f^{\ell,p}_\rho)$ is finite.  For some $\alpha < \delta$
we have $\eta \in {\gS}^{p_\alpha}_\ell$.  If $\rho \in
\cup\{\Lambda^{p_\beta}:\beta <  \delta\}$ then we apply $\odot_2$,
now if clause (a) there holds so $\alpha = \alpha(\rho) <
 \delta$ is well defined and we use $p_\alpha \in \bbQ$ and if
 clause (b) there holds then trivially Rang$(\eta) \subseteq
 u^{p_\delta} \subseteq \text { Dom}(f^{\ell,p_\delta}_\rho)$ 
 so assume $\rho \in \Lambda^{p_\delta}
\backslash \cup\{\Lambda^{p_\beta}:\beta < \delta\}$.

By $\odot_1$ we finish as in the case $\odot_2(b)$ holds.
\medskip

\noindent
\underline{Clause $(k)$}:

By the choice of ${\gS}^{p_\delta}$ and the proof of clause (g).
\medskip

\noindent
\underline{Checking $p_\alpha \le_{\bbQ} p_\delta$}: (where $\alpha < \delta$)

We should check that the pair $(p_\alpha,p_\delta)$ satisfies the
demands in $\circledast_2$ which is straight.
$\circledast_2$.

So we have proved $\circledast_4$. 
\mn
\begin{enumerate}
\item[$\circledast_5$]    if $\alpha < \mu$ then 
${\cI}^1_\alpha := \{p\in {\bbQ}:\alpha\in u^p\}$ is dense and
open as well as ${\cI}_* = \{p \in \bbQ$: if $\delta \in
u^p,\lambda|\delta$ and cf$(\delta) < \lambda$ then $\delta =
\sup(\delta \cap u)\}$.
\end{enumerate}
\mn
[Why?  Straight.  For the first, ${\cI}^1_\alpha$, given $p \in \bbQ$ we
define $q \in \bbQ$ by
\mn
\begin{enumerate}
\item[$(a)$]  $u^q$ is $u^p \cup\{\beta \le \alpha:\beta + \lambda =
\alpha + \lambda\}$, so clause $\circledast_1(a)$ holds
\sn
\item[$(b)$]  $<^q$ is the following linear order on $u^q$

$\alpha_1 < \alpha_2$ \Iff \, $\alpha_1 <^p \alpha_2$ or $\alpha_1 <
\alpha_2 \wedge \{\alpha_1,\alpha_2\} \nsubseteq
u^p \wedge \{\alpha_1,\alpha_2\} \subseteq u^q$
\sn
\item[$(c)$]  ${\gS}^q_\ell = {\gS}^p_\ell$ for $\ell=1,2$
\sn
\item[$(d)$]   $\Lambda^q = \Lambda^p$ and
\sn
\item[$(e)$]    $f^q_\rho = f^p_\rho$ for $\rho \in \Lambda^q$.
\end{enumerate}
\mn
Now check.

For the second, ${\cI}_*$ use the first and $\circledast_4$.]
\mn
\begin{enumerate}
\item[$\circledast_6$]    if $\varrho\in \Lambda^* :=
\{\rho:\rho$ is an increasing sequence of ordinals 
$<\lambda^+$ divisible by $\lambda$ 
of length $<\lambda\}$ then ${\cI}^2_\varrho
=\{p\in {\bbQ}:\varrho \in \Lambda^p\}$ is dense open.
\end{enumerate}
\mn
[Why?  Let $p\in {\bbQ}$, by $\circledast_5 + \circledast_4$ there is 
$q \ge p$ (from $\bbQ$) such that {\rm Rang}$(\varrho) \subseteq u^q$.
If $\varrho \in \Lambda^q$ we are done, otherwise define 
$q'$ as follows: 
$u^{q'}=u^q, <^{q'}=<^q, {\gS}^{q'}_\ell={\gS}^q_\ell, 
\Lambda^{q'}=\Lambda^q\cup \{\varrho \restriction \varepsilon:\varepsilon
\le \ell g(\varrho\}$  and if $i \le \ell g(\varrho),\varrho 
\restriction i\notin \Lambda^q$ then we let $f^{q'}_{\varrho \restriction i}
=\cup \{f^q_\rho:\rho\in \Lambda^q$ and $\rho 
\triangleleft \varrho \restriction i\}$.  We should check all the
clauses of $\circledast_1$ for ``$q \in \bbQ$" and e.g. clause (k)
there holds because $q$ satisfies clause $(\ell)$.  Then we should
check all the clauses of $\circledast_2$ for ``$q \le_{\Bbb Q} q'$"]
\mn
\begin{enumerate}
\item[$\circledast_7$]   if $\varrho$ is as in $\circledast_8$ and 
$\alpha<\lambda^+$ and $\ell\in \{1,2\}$ \then \,
\[
{\cI}^3_{\varrho,\alpha,\ell} = \{p\in {\bbQ}:\alpha\in \text{\rm Dom}
(f^{\ell,p}_\varrho) \text{ so } \varrho \in \Lambda^p, 
\alpha\in u^p\} \text{ is dense open}.
\]
\end{enumerate}
\mn
[Why?  By $\circledast_5 + \circledast_6$.]  
\mn
\begin{enumerate}
\item[$\circledast_8$]   if $p \in \bbQ$ and $\varrho \in
\Lambda^p$ \then \, for some $q$ we have $p \le_{\bbQ} q \wedge
f^q_\varrho \ne f^p_\varrho \wedge \{\alpha + \lambda:\alpha \in
u^q\} = \{\alpha + \lambda:\alpha \in u^p\}$.
\end{enumerate}
\mn
Why?  For each $\delta \in u \cap \text{ sup Rang}(\varrho)$ divisible by
$\lambda$ let $u_\delta = u \cap [\delta,\delta +\lambda)$.  So
$g_\delta := f^p_\rho \restriction u_\delta$ is a partial 
function from $u_\delta$
into $u_\delta$ and $f^p_\rho = \cup\{g_\delta:\delta$ as above$\}$.
Now, for $\delta$ as above we can find $f_\delta$ such that:
\mn
\begin{enumerate}
\item[$(a)$]  $f_\delta$ is a one-to-one function
\sn
\item[$(b)$]  $g_\delta = f^p_\varrho \restriction u_\delta \subseteq
f_\delta$
\sn
\item[$(c)$]  if $\alpha \in \text{ Dom}(f_\delta)$ iff $\alpha
\in u_\delta \vee f_\delta(\alpha) \in u_\delta$
\sn
\item[$(d)$]   Dom$(f_\delta) \backslash u_\delta$ is an initial segment
 $[\alpha^1_\delta,\alpha^2_\delta)$ of $[\delta,\delta + \lambda)
\backslash u_\delta$
\sn
\item[$(e)$]  Rang$(f_\delta) \backslash u$ is an initial segment
$[\alpha^2_\delta,\alpha^3_\delta)$ of $[\delta,\delta + \lambda)
\backslash u \backslash \text{ Dom}(f_\delta)$
\sn
\item[$(f)$]   $f_\delta$ maps $[\alpha^1_\delta,\alpha^2_\delta)$
onto $u_\delta \backslash \text{ Rang}(f^p_\varrho \restriction
u_\delta)$
\sn
\item[$(g)$]    $f_\delta$ maps $u_\delta \backslash 
\text{ Dom}(f^p_\varrho \restriction u_1)$ onto
$[\alpha^2_\delta,\alpha^3_\delta)$.
\end{enumerate}
\mn
Now we can find a linear order $<_1$ on $u_\delta \cup
[\alpha^1_\delta,\alpha^3_\delta]$ such that $f_\delta$ is order preserving
(as the class of linear orders has amalgamation).

Lastly, we define $q$:
\mn
\begin{enumerate}
\item[$(\alpha)$]   $u^q = u^p \cup
\{[\alpha^1_\delta,\alpha^3_\delta):\delta$ as above$\}$
\sn
\item[$(\beta)$]   $<^q$ is defined by $\alpha <^q \beta$ iff
$(\exists \delta)(\alpha <_\delta \beta)$ or $\alpha + \lambda \le \beta$
\sn
\item[$(\gamma)$]  $\Lambda^q = \Lambda^p$
\sn
\item[$(\delta)$]   ${\gS}^q_\ell = {\gS}^p_\ell \cup
\{\langle f^{3-\ell}_\rho(\eta(n)):n < \omega\rangle:\ell
\in\{1,2\},\rho \in \Lambda^p$ and $\eta \in {\gS}^p_{3-\ell}\}$.
\end{enumerate}
\mn
Now we have to check $q \in \bbQ$, i.e. all the clauses of
$\circledast_1$.  This is straight; e.g. for clause (c), assume $\eta
\ne \nu \in {\gS}^q_\ell$ and we have to prove that Rang$(\eta)
\cap \text{ Rang}(\nu)$ is finite.

Now we have four cases: first $\eta,\nu \in {\gS}^p_\ell$, so use $p
\in \bbQ$, clause $\circledast_1(c)$ for $\ell$.  
Second, $\eta,\nu \in {\gS}^q_\ell
\backslash {\gS}^p_\ell$, so $\eta,\nu$ are images by
$f^{3-\ell,q}_\rho$ of members of ${\gS}^p_{3-\ell}$, as this
function is one-to-one, this follows from $p,{\gS}^p_{3-\ell}$
satisfying clause $\circledast_1(c)$.  Third, $\eta \in {\gS}^p_\ell
\wedge \nu \in {\gS}^q_\ell \backslash {\gS}^p_\ell$, then
$\nu = \langle f^{3-\ell,q}_\rho(\nu'(n)):n < \omega\rangle$ for some
$\nu' \in {\gS}^p_{3-\ell}$ satisfying Rang$(\nu') \nsubseteq
\text{ Dom}(f^{3-\ell,p}_\rho)$, hence for some $n_* < \omega$ we have
$n \in [n_*,\omega) \Rightarrow \nu'(n) \notin \text{ Dom}(f^{3-\ell,p}_\rho)
\Rightarrow \nu(n) \notin u^p$ but Rang$(\eta) \subseteq u^p$ so we
are done.  Fourth, $\eta \in {\gS}^q_\ell \backslash 
{\gS}^p_\ell \wedge \nu \in {\gS}^p_\ell$ the proof is dual.

The proof of clause (g) is similar.

Also we have to check that $p \le_{\bbQ} q$, i.e. all the clauses of
$\circledast_2$ for the pair $(p,q)$.  This is straight, clause (f) is
sproved as in the proof of $\circledast_1(c)$ above and clause (h)
holds by our choice of the $f_\delta$'s.

Now check that $q$ is as required.]

Let
\mn
\begin{enumerate}
\item[$\oplus_1$]   $\bbQ^+ = \{p \in \bbQ$: if $\ell \in
\{1,2\}$ and $\rho \in \Lambda^p$ then Dom$(f^{\ell,p}_\rho) = u^p \cap
\text{ sup Rang}(\rho)\}$
\sn
\item[$\oplus_2$]   $\bbQ^+$ is a dense subset of $\bbQ$, 
moreover $(\forall p \in \bbQ)(\exists q \in \bbQ^+)(p \le q
\wedge \{\alpha + \lambda:\alpha \in u^q\} = \{\alpha +
\lambda:\alpha \in u^p\}]$.
\end{enumerate}
\mn
[Why?  Let $p \in \bbQ,\kappa = |\Lambda|,\delta = \kappa \times
\kappa$ and $\{\rho_i:i < i_* < \lambda\}$ list $\Lambda^p$ each
appearing unboundedly often.  We choose $p_i$ by induction on $i \le
\delta$ such that
\mn
\begin{enumerate}
\item[$(a)$]   $p_i \in \bbQ$
\sn
\item[$(b)$]   $j < i \Rightarrow p_i \le_Q p_j$
\sn
\item[$(c)$]   $p_0 = p$
\sn
\item[$(d)$]  $\Lambda^{p_i} = \Lambda^p$
\sn
\item[$(e)$]  $f^{p_{i+1}}_{\rho_i} \ne f^{p_i}_{\rho_i}$
\sn
\item[$(f)$]   $\{\alpha + \lambda:\alpha \in u^{p_i}\} =
\{\alpha + \lambda:\alpha \in u^p\}$.
\end{enumerate}
\mn
For $i=0$ use clause (c) for $i$ limit use $\circledast_4$, for
$i=j+1$ use $\circledast_8$.  Now $p_\delta$ is as required.]
\mn
\begin{enumerate}
\item[$\oplus_3$]    for $p \in \bbQ$ and $\delta <
\lambda^+$ divisible by $\lambda,p \restriction \delta$ is naturally
defined, belongs to $\bbQ$ and $u^p \subseteq \delta \Rightarrow p
\rest \delta =p$ and $p \restriction \delta \le_{\bbQ} p$, where 
$q = p \restriction \delta$ be defined by:
\begin{enumerate}
\item[$(a)$]  $u^q = u^p \cap \delta$ 
\sn
\item[$(b)$]   $<^q = <^p \restriction \delta$
\sn
\item[$(c)$]   ${\gS}^q_\ell = \{\eta \in {\gS}^p_\ell$:Rang$(\eta) 
\subseteq \delta\}$
\sn
\item[$(d)$]    $\Lambda^q = \{\rho \in \Lambda^p:
\text{\rm sup Rang}(\rho) \le \delta\}$
\sn
\item[$(e)$]   $\bar f^q = \langle f^q_\rho:\rho \in
\Lambda^q\rangle$ where $f^q_\rho = f^p_\rho$.
\end{enumerate}
\end{enumerate}
\mn
[Why?  Check.]
\mn
\begin{enumerate}
\item[$\oplus_4$]  if $\delta < \lambda^+$ is divisible by
$\lambda,p \in \bbQ^+$ and $(p \rest \delta) \le_{\bbQ} q \in 
\bbQ^+$ but $u^q \subseteq \delta$ \then \, $p,q$ are compatible in
$\bbQ$, moreover has a common upper bound $r = p+q$ such that $r \rest
\delta = q \wedge u^r = u^p \cup u^q$.
\end{enumerate}
\mn
[Why?  Note that if $\rho \in \Lambda^p \cap \Lambda^q$ then sup
Rang$(\rho) \le \delta$ by clause (i) + (j) of $\circledast_1$; also
$\Lambda^p \cap \Lambda^q =\Lambda_{p \rest \delta}$.  We
define $r$ as follows:
\mn
\begin{enumerate}
\item[$(a)$]   $u^r = u^p \cup u^q$
\sn
\item[$(b)$]   $\le^r$ is defined by: for $\alpha,\beta \in u^r$ we
have $\alpha <^2 \beta$ iff $\alpha + \lambda \le \beta$ or $\alpha
<^q \beta$ or $\alpha <^p \beta$
\sn
\item[$(c)$]   ${\gS}^r_\ell$ is ${\gS}^p_\ell \cup {\gS}^q_\ell$ 
for $\ell=1,2$
\sn
\item[$(d)$]  $\Lambda^r = \Lambda^p \cup \Lambda^q$
\sn
\item[$(e)$]   $\bar f^r = \langle f^r_\rho:\rho \in
\Lambda^r\rangle$ where $f^r_\rho$ is:
\begin{enumerate}
\item[$\bullet$]   $f^q_\rho$ \when \, $\rho \in \Lambda^q$
\sn
\item[$\bullet$]   $f^p_\rho \cup \bigcup\{f^q_{\rho \rest i}:i \le \ell
g(\rho)$ and $\rho \rest i \in \Lambda^q\}$ \when \, $\rho \in
\Lambda^p \backslash \Lambda^q$.
\end{enumerate}
\end{enumerate}
\mn
Why $r \in \bbQ$?  We should check all the clauses in
$\circledast_1$, which are easy.  E.g. in clause (c), $\eta \ne \nu
\in {\gS}^r_\ell \Rightarrow \aleph_0 > |\text{Rang}(\eta) \cap 
\text{ Rang}(\nu)|$, the only new case is $\eta \in {\gS}^p_\ell
\Leftrightarrow \nu \notin {\gS}^p_\ell$ so \wilog \, $\eta \in
{\frak S}^p_\ell \backslash {\gS}^q_\ell \wedge \nu \in 
\{\gS\}^q_\ell$, hence sup$(\eta) >
\delta$ hence Rang$(\eta) \cap \delta$ is finite but Rang$(\nu)
\subseteq u^q \subseteq \delta$.

Also clauses (g) + (k) should be checked only when $f^r_\rho$ is new
so necessarily $\rho \in \Lambda^p$ so $f^r_\rho = f^p_\rho \cup
\bigcup \{f^q_{\rho \rest i}:\rho \rest i \in \Lambda^q\}$, but
recalling that any $\eta \in {\gS}^r_\ell$ is an increasing
$\omega$-sequence, clearly if sup Rang$(\eta) > \delta$ we use ``$p$
satisfies clauses (g) + (k)" and if sup Rang$(\eta) \le \delta$ we use
``$q$ satisfies clauses (g) +(k) and $(\ell)$".

Why $p \le_{\bbQ} r \wedge p \le_{\bbQ} r$?  We should check all
the clauses in $\circledast_2$ for both pairs.  They are easy,
e.g. clause (f)  holds because: if $\eta \in {\gS}^r_\ell
\backslash {\gS}^q_\ell$ then $\eta \in {\gS}^p_\ell
\backslash S^q_\ell$ hence sup Rang$(\eta) > \delta$ and it should be
clear; if $\eta \in {\gS}^r_\ell \backslash {\gS}^p_\ell$ then
$\eta \in{\gS}^q_\ell \backslash {\gS}^p_\ell$ and we can use
$p \rest \delta \le_{\bbQ} q$, i.e. clause (f) for this pair.

Concerning clause (g) for $p \le_{\bbQ} r$, recall that $p,q \in
\bbQ^+$ so $\ell \in \{1,2\} \wedge \rho \in \Lambda^p \Rightarrow u^p =
\text{ Dom}(f^{\ell,p}_\rho) \subseteq \text{ Dom}(f^{\ell,r}_\rho)$
so clause (g) is O.K. and similarly clause (g) for $q \le_{\bbQ} r$.]
\mn
\begin{enumerate}
\item[$\oplus_5$]  $\bbQ$ satisfies the $\lambda^+$-c.c.
\end{enumerate}
\mn
[Why?  Let $p_\alpha \in \bbQ$ for $\alpha < \lambda^+$, so by
$\circledast_{10}$ there are $q_\alpha$ such that $p_\alpha \le_{\bbQ} 
q_\alpha \in \bbQ^+$, now use the $\Delta$-sytem lemma 
that is first $S^{\lambda^+}_\lambda = \{\delta < \lambda^+$:
cf$(\delta) = \lambda\}$; now $\delta \in S^{\lambda^+}_\lambda
\Rightarrow p \rest \delta \in \bbQ \wedge \sup(u^{p \rest \delta}) <
\delta$ and $\lambda \ge |\{p \in \bbQ:u^p=u\}$ for any $u$.  Hence
for some stationary $S \subseteq S^{\lambda^+}_\lambda$ and $p_*$ we
have $\delta \in S \Rightarrow q_\delta \rest \delta = p_*$ and
$\delta_1 < \delta_2 \in S \Rightarrow \sup(u^{q_{\delta_2}}) <
\delta_2$.  So for any $\delta_2 < \delta_2$ from $S$ by $\oplus_4$
the condition $q_{\delta_1},q_{\delta_2}$ are compatible.]
\mn
\begin{enumerate}
\item[$\boxplus_1$]  define 
$\name{\bold J}_\ell \in K^{\text{oi}}_\mu$, a ${\bbQ}$-name as follows:
\begin{enumerate}
\item[$(a)$]   $Q^{\name{\bold J}_\ell} =\mu$
\sn
\item[$(b)$]   ${\gS}^{\name{\bold J}_\ell} =
\cup \{{\gS}^p_\ell:p\in \name G_{\Bbb Q}\}$
\sn
\item[$(c)$]   $<^{\name{\bold J}_\ell} =
\cup \{ <^p:p\in \name G_{\Bbb Q}\}$
\sn
\item[$(d)$]   $F^{\name{\bold J}_\ell}_n$ is a 
unary function, the identity on $\lambda^+$ and 
\sn
\item[$(e)$]   $\eta \in {\gS}^{\name{\bold J}_\ell} 
\Rightarrow F^{\name{\bold J}_n} (\eta)=\eta(n)$
\end{enumerate}
\item[$\boxplus_2$]   for $\ell \in \{1,2\}$ and
$p\in {\bbQ}$ let ${\bold J}^p_\ell \in K_{\text{oi}}$ be defined as follows:
\begin{enumerate}
\item[$(a)$]   $\bold J^p_\ell$ has universe $u^p \cup {\gS}^p_\ell$
\sn
\item[$(b)$]  $<^{{\bold J}_\ell} = <^p$
\sn
\item[$(c)$]   $Q^{{\bold J}^p_\ell} = u^p$
\sn
\item[$(d)$]  $F^{{\bold J}^p_\ell}_n (\eta)= \eta(n)$
\end{enumerate}
\item[$\boxplus_3$]    $(a) \quad \Vdash_{\bbQ}$ 
``$\name{\bold J}_\ell \in K^{\text{oi}}_{\lambda^+}$"
\sn
\item[${{}}$]   $(b) \quad \Vdash_{\bbQ}$ ``for each $\delta <
\lambda^+$ divisible by $\lambda$ the linear order $([\delta,\delta +
\lambda),<^{\name{\bold J}_\ell} \rest (\delta,\delta + \lambda))$ is a
saturated linear order and $\alpha + \lambda \le \beta < \lambda^+
\Rightarrow \alpha <^{\name{\bold J}_\ell} \beta$"
\sn
\item[${{}}$]   $(c) \quad p \in \bbQ \Rightarrow p
\Vdash_{\bbQ} ``\bold J^p_\ell \subseteq \name{\bold J}_\ell$ 
for $\ell=1,2"$.
\end{enumerate}
\mn
[Why?  Think]
\mn
\begin{enumerate}
\item[$\boxplus_4$]    if $\delta < \lambda^+$ is divisible by
$\lambda$ then $\Vdash ``\name{\bold J}_\ell \restriction \delta \in
K^{\text{oi}}_\lambda$ where
${\bold J}_\ell \restriction \delta= ((\delta\cup (P^{{\bold J}_\ell}
\cap {}^\omega\delta),Q^{\bold J_\ell} \cap \delta,P^{\bold J_\ell} 
\restriction \delta,F_n^{{\bold J}_\ell} \restriction (\delta \cup 
(P^{{\bold J}_\ell}\cap {}^{\omega}\delta)))_{n < \omega}"$
\sn
\item[$\boxplus_5$]   $\Vdash_{\bbQ}
``\text{\rm EM}_{\tau(T)} (\name{\bold J}_1,\Phi),
\text{\rm EM}_{\tau(T)} (\name{\bold J}_2,\Phi)$ 
are {\rm EF}$^+_{\lambda,\lambda^+}$-equivalent
\end{enumerate}
\mn
(so the games of length $<\lambda$, and the player INC chooses sets of  
cardinality $<\lambda^+$).

\noindent
[Why?  To show the EF$^+_{\lambda,\lambda^+}$-equivalence, it suffices
to show that $\Vdash_{\bbQ} ``\name{\bold J}_1,\bold J_2$ are
EF$_{\lambda,\lambda^+}$-equivalent" by \ref{1a.24} as $\lambda \ge
\aleph_1 + |T_1|$.  From $\circledast_6$, recall 
$\Lambda^*=\{\rho:\rho$ is an increasing sequence
of ordinals $<\lambda^+$ divisible by $\lambda$ of length 
$<\lambda\}$, (is the same in ${\bold V}$ and ${\bold V}^{\bbQ}$).
For $\rho\in \Lambda^*$ let $\name f_\rho =
\cup\{f^p_\rho:\rho \in \name G,p \in \Lambda^p\}$
and by $\circledast_1(f) + (j) \and \circledast_2(e)$ easily 
$\Vdash_{\bbQ} ``\name f_\rho$ a partial isomorphism 
from $\name{\bold J}_1 \restriction \text{\rm sup Rang} 
(\rho)$ into $\name{\bold J}_2 \restriction 
\text{\rm sup Rang}(\rho)"$, see Definition inside $\boxplus_4$.

Now $\Vdash_{\bbQ} ``\text{\rm Dom}(\name f_\rho) =
\text{\rm sup Rang}(\rho)"$ as if $\name G \subseteq
\bbQ$ is generic over $\bold V$, for any $\alpha < \text{\rm sup
Rang}(\rho)$ for some $p \in \bold G$ we have $\alpha \in u^p \wedge
\rho \in \Lambda^p$ by $\circledast_7$ and there is $q$ such that 
$p \le q \in \bold G,p \ne q$ by $\circledast_8$, 
so recalling $\circledast_2(g)$ we are
done.

Similarly $\Vdash_{\bbQ} ``\text{\rm Rang}(\name f_\rho) 
= \text{\rm sup Rang}(\rho)"$.  

Also $\rho \triangleleft \varrho \Rightarrow \Vdash_{\bbQ} 
\name f_\rho \subseteq \name f_\varrho$.
For the EF$^+$-version we have to analyze dependence relations, which
is straight as in the proof in \ref{4d.7}.
So $\langle f_\rho:\rho\in \Lambda^*\rangle$ exemplify the equivalence.] 
\mn
\begin{enumerate}
\item[$\boxplus_6$]   $\Vdash_{\bbQ} ``\name M_1 = 
\text{ EM}_{\tau(T)}(\name{\bold J}_1,\Phi), 
\name M_2 = \text{ EM}_{\tau(T)} 
(\name{\bold J}_2,\Phi)$ are not isomorphic".
\end{enumerate}
\mn
Why?  Let $\name M^+_\ell = \text{ EM} 
(\name{\bold J}_\ell,\Phi)$ so $\name M^+_\ell
\rest \tau(T) = \name M_\ell$ for $\ell=1,2$, and assume
toward contradiction that $p\in {\bbQ}$, and 
$p \Vdash_{\bbQ}$ ``$\name g$ is an isomorphism from $\name M_1$ 
onto $\name M_2$".
For each $\delta\in S^{\lambda^+}_\lambda:=
\{\delta<\lambda^+:\text{cf}(\delta)=\lambda\}$ by $\circledast_4$ we can find 
$p_\delta\in {\bbQ}$ above $p$ and $g_\delta$ such that: 
\mn
\begin{enumerate}
\item[$\boxdot_1$]  $(a) \quad p \le p_\delta,\delta\in u^{p_\delta}$
\sn
\item[${{}}$]  $(b) \quad p_\delta \Vdash 
``g_\delta$ is $\name g \restriction 
\text{\rm EM}({\bold J}^{p_\delta}_1,\Phi)"$
\sn
\item[${{}}$]  $(c) \quad g_\delta$ is an isomorphism from
EM$_{\tau(T)} ({\bold J}^p_1,\Phi)$ onto EM$_{\tau(T)}({\bold J}^p_2,\Phi)$. 
\end{enumerate}
\mn
We can find stationary $S \subseteq S^{\lambda^+}_\lambda$ and $p^{*}$
such that
\mn
\begin{enumerate}
\item[$\boxdot_2$]   $(a) \quad p_\delta \restriction \delta$, 
 defined in $\oplus_3$ is $p^*$ for $\delta \in S$
\sn
\item[${{}}$]   $(b) \quad$ for $\delta_1,\delta_2\in S$, 
 $u^{p_{\delta_1}},u^{p_{\delta_2}}$ has the same order type and 
the order

\hskip25pt  preserving mapping $\pi_{\delta_1,\delta_2}$ 
from $u^{p_{\delta_2}}$ onto $u^{p_{\delta_1}}$ induce an

\hskip25pt  isomorphism from $p_{\delta_2}$ onto $p_{\delta_1}$
\sn
\item[${{}}$]  $(c) \quad$ if $\delta_1 < \delta_2 \in S$ then
sup$(u^{p_{\delta_1}}) < \delta_2$.
\end{enumerate}
\mn
Now choose $\eta^*= \langle \delta^*_n:n<\omega\rangle$ such that
\mn
\begin{enumerate}
\item[$\boxdot_3$]  $(a) \quad \delta^*_n < \delta^*_{n+1}$
\sn
\item[${{}}$]  $(b) \quad \delta^*_n =\text{\rm sup}(S\cap
\delta^*_n)$ and $\delta^*_n \in S$
\sn
\item[${{}}$]  $(c) \quad$ let $\delta^* = 
\text{ sup}\{\delta^*_n:n < \omega\}$.
\end{enumerate}
\mn
We define $q \in {\bbQ}$ as follows
\mn
\begin{enumerate}
\item[$\boxdot_4$]   $(a) \quad u^q =
\cup\{p_{\delta^*_n}:n<\omega\}$
\sn
\item[${{}}$]  $(b) \quad <^q= \{(\alpha,\beta):\alpha <^{p_{\delta^*_n}}
\beta$ for some $n$ or $\alpha + \lambda \le \beta \wedge
\{\alpha,\beta\} \subseteq u^q$, 

\hskip25pt equivalently for some $m < n,
\alpha\in u^{p_{\delta^*_m}} \setminus \delta^*_m$ and $\beta\in 
u^{p_{\delta^*_n}} \setminus \delta^*_n\}$
\sn
\item[${{}}$]   $(c) \quad {\gS}^q_1 = 
\cup\{{\gS}_1^{p_{\delta^*_n}}:n<\omega\}\cup \{\eta^*\}$
\sn
\item[${{}}$]   $(d) \quad {\gS}^q_2 = 
\cup\{{\gS}_2^{p_{\delta^*_n}}:n<\omega\}$
\sn
\item[${{}}$]   $(e) \quad \Lambda^q =
\cup \{ \Lambda^{p_{\delta^*_n}}:n < \omega\}$
\sn
\item[${{}}$]  $(f) \quad f^q_\rho = f^{p_{\delta^*_n}}_\rho$ if 
$\rho\in \Lambda^{p_{\delta^*_n}}$.
\end{enumerate}
\mn
So there is a pair $(q_*,g^+)$ such that:
\mn
\begin{enumerate}
\item[$\boxdot_5$]   $(a) \quad q \le_{\bbQ} q_*$
\sn
\item[${{}}$]  $(b) \quad q_* \Vdash_{\bbQ} ``g^+ = \name g 
\restriction \text{ EM}(\bold J^{q_*}_1,\Phi)$
\sn
\item[${{}}$]   $(c) \quad g^+$ is an isormorphism from
EM$_{\tau(T)}(\bold J^{q_*}_1,\Phi)$ onto EM$_{\tau(T)}(\bold
J^{q_*}_2,\Phi)$.
\end{enumerate}
\mn
So $g^+(a_{\eta^*}) \in \text{ EM}(\bold J^{q_*}_2,\Phi)$ hence is of the form
$\sigma^{M^+_2}(a_{t_0},\dotsc,a_{t_{n-1}})$ for some
$t_0,\dotsc,t_{n-1} \in \bold J^{q_*}_2$ and a $\tau_\Phi$-term 
$\sigma(x_0,\dotsc,x_{n-1})$. 

Note that by the definition of $\le_{\bbQ}$ in $\circledast_2$:
\mn
\begin{enumerate}
\item[$\boxdot_6$]   if $\eta \in {\gS}^{q_*}_2$ then Rang$(\eta)
\cap u^q$ is bounded in $\delta^*$.
\end{enumerate}
\mn
[Why?  If $\eta \in {\gS}^q_2$ this holds by our choice of $q$ and
if $\eta \in {\gS}^{q_*}_2 \backslash {\gS}^q_2$ then
Rang$(\eta) \cap u^q$ is finite so as $u^q \subseteq \delta$ it
follows that Rang$(\eta) \cap u^q$ is bounded in $\delta^*$.]

We can find $n(*) < \omega$ such that:
\mn
\begin{enumerate}
\item[$\boxdot_7$]  for each $k < n$ and $\ell < n$ we have
\begin{enumerate}
\item[$(a)$]  if $t_\ell \in Q^{\bold J^{q_*}_2}$, i.e. $t_\ell
\in u^{q_*} \subseteq \lambda^+$ then $t_\ell \le^q \delta^*_{n(*)}$
or $\delta^* \le t_\ell$ (hence $\bigwedge\limits_{n} \delta^*_n \le^q t_\ell$)
\sn
\item[$(b)$]  if $t_\ell \in P^{\bold J^{q_*}_2}$, i.e. $t_\ell
\in {\frak S}^{q_*}_2$ then $\{F^{\bold J^{q_*}_2}_n(t_\ell):n <
\omega\}$ is disjoint to $[\delta^*_{n(*)},\delta^*) \cap u^q$.
\end{enumerate}
\end{enumerate}
\mn
Now using ``$T$ is stable", the rest is as in \ref{2b.4}, 
\ref{2b.6}.
\end{PROOF}

\begin{discussion}
\label{3c.12}
(2012.11.23)
1) Can we do it in ZFC?  It is natural to use $\langle W_\alpha:\alpha
\in S^\mu_{\aleph_0}\rangle$ be stationary pairwise almost disjoint,
see \cite{Sh:F980}.

\noindent
2) Instead of ``not strongly stable" it suffices to assume ``not
   strongly$^2$ stable", see \cite{Sh:863}.  In \cite{Sh:F930} even
   much less.
\end{discussion}
\newpage

\section {Theories with order} 

Recall from \cite[3.19]{HySh:529}:
\begin{claim}
\label{4d.2}
If $\lambda = \mu^+$, {\rm cf}$(\mu),\lambda
= \lambda^{< \kappa},\kappa = \text{\rm cf}(\kappa) < \kappa(T)$ and
$T$ is unstable \then \, there are 
{\rm EF}$_{\mu \times \kappa,\lambda^+}$-equivalent
non-isomorphic models of $T$ of cardinality $\lambda$.
\end{claim}

\noindent
The new point in \ref{4d.7} is the EF$^+$ rather than EF.
\begin{claim}
\label{4d.7}
Assume $\lambda = \lambda^{< \theta}$ and
$\lambda$ is regular uncountable, $T
\subseteq T_1$ are complete first order theories of cardinality $< \lambda$.

\noindent
1) If $T$ is unstable \then \, there are models $M_1,M_2$ 
of $T_1$ of cardinality $\lambda^+$, 
{\rm EF}$^+_{\lambda,\theta,\lambda^+}$-equivalent
 with non-isomorphic $\tau_T$-reducts. 

\noindent
2) Assume $\Phi \in \Upsilon^{\text{or}}_\kappa$ 
is proper for linear orders, $\bar \sigma = \langle
\sigma_i(x):i <i(*)\rangle$ a sequence of terms from $\tau_\Phi,\bar x^\ell =
\langle x^\ell_i:i <i(*)\rangle,i(*) < \lambda,\varphi(\bar x^1,
\bar x^2)$ is a formula in $\bbL(\tau_T),\tau \le \tau_T$ (any
logic) and for every linear order $I$ letting $M = 
\text{\rm EM}(I,\Phi),\bar b_t = \langle \sigma^M_i(a_t):i<i(*)\rangle$ we have
$(M \restriction \tau) \models \varphi[\bar b_s,\bar
b_t]^{\text{if}(s<t)}$ for every $s,t \in I$.  \Then \, there are
linear orders $I_1,I_2$ of cardinality $\lambda^+$ such that $M_1,M_2$
are {\rm EF}$^+_{\lambda,\theta,\lambda^+}$-equivalent but not isomorphic where
$M_\ell = \text{\rm EM}_\tau(I_\ell,\Phi)$ for $\ell=1,2$.

\noindent
3) If every {\rm EM}$_\tau(I,\Phi)$ is a model of $T_1$ \then \, in (2) the
models $M_1,M_2$ are in {\rm PC}$(T_1,T)$.
\end{claim}

\begin{PROOF}{\ref{4d.7}}
1) Let $\varphi(\bar x,\bar y) \in \bbL(\tau_T)$ order
some infinite subset of ${}^m M$ for some $M \models T$.  

Let $\Phi$ be as in Definition \ref{2b.2A}, i.e.
\cite[Ch.VII,VIII]{Sh:c}, i.e. proper for linear
orders such that $\tau_{T_1} \subseteq \tau(\Phi),|\tau(\Phi)| =
|T_1|$ and for every linear order $I$,EM$(I,\Phi)$ (we allow the
skeleton to consist of $m$-tuples rather than elements) is a model of $T_1$
satisfying $\varphi[\bar a_s,\bar a_t]$ iff $s <_I t$.  Now we can
apply part (2) with $i(*) = m$.

\noindent
2) We choose $I$ such that
\mn
\begin{enumerate}
\item[$\circledast$]  $(a) \quad I$ is a linear order of
cardinality $\lambda$ (yes, not $\lambda^+$)
\sn
\item[${{}}$]  $(b) \quad$ if $\alpha,\beta \in (1,\lambda]$ then 
$(I \times \alpha) + (I \times \beta)^* \cong I$ 
(equivalently every $\alpha$,

\hskip25pt $\beta \in [1,\lambda^+))$
\sn
\item[${{}}$]  $(c) \quad I$ is isomorphic to its inverse
\sn
\item[${{}}$]  $(d) \quad I$ has cofinality $\lambda$.
\end{enumerate}
\mn
For every $S \subseteq S^{\lambda^+}_\lambda = \{\delta <
\lambda^+:\text{cf}(\delta) = \lambda\}$ we define $I_S = 
\sum\limits_{\alpha < \lambda^+} I_{S,\alpha}$ where $I_{S,\alpha}$ is
isomorphic to $I$ if $\alpha \in \lambda^+ \backslash S$ and isomorphic
to the inverse of $I \times \omega$ otherwise.  Now
\mn
\begin{enumerate}
\item[$\circledast_2$]  if $S_1,S_2 \subseteq
S^{\lambda^+}_\lambda$ \then \, the models 
EM$(I_{S_2},\Phi)$, EM$(I_{S_1},\Phi)$ are
EF$^+_{\lambda,\theta,\lambda^+}$-equivalent.
\end{enumerate}
\mn
[Why?  Let $J_{\ell,\gamma} = \sum\limits_{\alpha < \gamma}
I_{S_\ell,\alpha}$.  Let ${\cF} := \{f$: for some non-zero ordinal
$\gamma < \lambda^+,f \in {\cF}_\gamma$ and $[\gamma \in S_1
\Leftrightarrow \gamma \in S_2]\}$ where ${\cF}_\gamma := \{f$ is an
isomorphism from $\sum\limits_{\alpha < \gamma} I_{S_1,\alpha}$ onto
$\sum\limits_{\alpha < \gamma} I_{S_2,\alpha}\}$.

Now
\mn
\begin{enumerate}
\item[$(*)_1$]    ${\cF}_\gamma \ne \emptyset$ for $\gamma <
\lambda^+$
\sn
\item[$(*)_2$]    if $f \in {\cF}_\gamma$ and $[\gamma \in S_1
\equiv \gamma \in S_2]$ and $\gamma < \beta < \lambda$ then $f$ can be
extended to some $g \in {\cF}_\beta\}$
\sn
\item[$(*)_3$]   if $\gamma < \lambda,X_\ell \subseteq I_\ell$ has
cardinality $< \lambda^+$ for $\ell=1,2$ \then \, for some successor
$\beta,\gamma < \beta < \lambda^+$ and $X_\ell \subseteq
J_{\ell,\beta}$ for $\ell=1,2$
\sn
\item[$(*)_4$]  if $\gamma_i \in S_1 \Leftrightarrow \gamma_i \in
S_2$ for $i < \delta,\delta$ a limit ordinal $< \lambda$ and $\langle
\gamma_i:i < \delta\rangle$ is increasing \then \, $\gamma_\delta :=
\cup\{\gamma_i:i < \delta\}$ satisfies $\gamma_\delta \in S_1 \equiv
\gamma_\delta \in S_1$.
\end{enumerate}
\mn
Lastly, we have to deal with case 2 in Definition \ref{1a.9}(2) so
assume
\mn
\begin{enumerate}
\item[$(*)_5$]   $f_* \in {\cF}_{\gamma_*},[\gamma_* \in S_1 \equiv
\gamma_* \in S_2]$ and $\bold R_\ell \subseteq {}^{\theta >}(M_\ell)$
for $\ell=1,2$ are as there for $f_*$ 
\end{enumerate}
\mn
This holds because the strategy is simple, e.g. with no memory.  Now
if $f$ does not map the definition of $\bold R_1$ in $M_1$ to the
definition of $\bold R_2$ in $M_2$ we can use subcase 2B there, so we
assume this does not occur.  Let $\ell \in \{1,2\}$.
\mn
\begin{enumerate}
\item[$(*)_6$]  Let $\bold e_\ell = 
\{(\bar s,\bar t):\bar s,\bar t \in {}^{\theta >}(I_\ell)$ 
and some automorphism of $I_\ell$ over $I_{\ell,\gamma_*}$ maps $\bar s$
to $\bar t\}$.
\sn
\item[$(*)_7$]   Let $Y_\ell$ be the set of $\bold e_\ell$-equivalence 
classes.
\end{enumerate}
\mn
Note  
\mn
\begin{enumerate}
\item[$\odot_1$]   for $\ell \in \{1,2\},n< \omega$ and 
$\bold y_0,\dotsc,\bold y_n \in Y_\ell$ the following are equivalent:
\begin{enumerate}
\item[$(a)$]    some $\bar a \in \bold y_n$ depend (by $\bold R_1$)
on $\bold y_0 \cup \ldots \cup \bold y_{n-1}$
\sn
\item[$(b)$]  every $\bar a \in \bold y_n$ depend (by $\bold R_1$) on  
$\bold y_0 \cup \ldots \cup \bold y_{n-1}$.
\end{enumerate}
\end{enumerate}
\mn
So $\bold R_1$ induce a 1-dependence relation on $Y_1$, so let $\langle
\bold y_i:i < i(*)\rangle$ be a maximal independent subset of $Y_1$
such that $[i <
i(*) \wedge \bar a \in \bold y_i \Rightarrow \bar a$ does not depend
on $\cup\{\bold y_j:j < i(*),j \ne i\}$.

So
\mn
\begin{enumerate}
\item[$\odot_2$]   it is enough to deal with one $\bold y_i$.
\end{enumerate}
\mn
Now we can find $\bar t_{i,\gamma} \in \bold y_i$ such that
Rang$(\bar t_{i,\gamma}) \backslash I_{\ell,\gamma_*} 
\subseteq I_{\ell,\gamma +2} \backslash I_{\ell,\gamma +1}$ for each
$\gamma \in [\gamma_*,\lambda^+)$ as $I_1$ has enough automorphisms
\mn
\begin{enumerate}
\item[$\odot_3$]   if $\{\bar t_{i,\gamma}:\gamma \in
[\gamma_*,\lambda^+)\}$ is not $\bold R_1$-independent, then
dim$(\bold y_i)$ is finite, in fact 1 or 0.
\end{enumerate}
\mn
So we choose $\beta_*$ such that
\mn
\begin{enumerate}
\item[$\odot_4$]    $\gamma_* < \beta_* < \lambda^+$ and $\beta_*
\in S_1 \equiv \beta_* \in S_2$ and for every $i < i(*)$, if
dim$(X_{\bold y_i})$ is finite then $\bold y_i$ has a maximal $\bold
R_1$-independent set included in ${}^{\varepsilon(\bold
y_i)}(J_{1,\beta_*})$.
\end{enumerate}
\mn
[Why possible?  Because for any such $\beta_*$ is an automorphism of
$I_2$ over $J_{1,\gamma_*}$ mapping $I_{\beta_* +2}$ onto $I_{\gamma_*
+2}$.]

Let $g \in {\cF}_{\beta_*}$ extend $f$ and using it we can
choose $\langle (\bar a^1_\zeta,\bar a^2_\zeta):\zeta <
\zeta^*\rangle$ as required. 
\mn
\begin{enumerate}
\item[$\circledast_3$]   if $S_1,S_2 \subseteq 
S^{\lambda^+}_\lambda$ and $S_1 \backslash S_2$ is stationary,
then EM$_\tau(I_{S_1},\Phi)$, EM$_\tau(I_{S_2},\Phi)$ are not isomorphic.
\end{enumerate}
\mn
[Why?  By the proof in \cite[Ch.III,\S3]{Sh:300} (or \cite[Ch.III,\S3]{Sh:e} =
\cite[\S3]{Sh:E59}) only easier.  
In fact, immitating it we can represent the invariants from there.  If
$k=2$ we have to work somewhat more.]

\noindent
2) As in \cite{HySh:529}.

\noindent
3) Obvious. 
\end{PROOF}

\begin{conclusion}
\label{4d.14}
Assume $T$ is a (first order complete) theory.

\noindent
1) If $T$ is unstable, \then \, $(T,*)$ is fat. 

\noindent
2)  If $T$ is unstable or stable with DOP, or stable with OTOP, 
\then \, $T$ is fat.

\noindent
3) For every $\mu$ there is a $\mu$-complete, class forcing $\bbP$
such that in $\bold V^{\bbP}$ we have: if $T$ is not strongly dependent
or just not strongly stable \then \, $T$ is fat, moreover $(T,*)$ is fat. 
\end{conclusion}

\begin{PROOF}{\ref{4d.14}}
 1) By \ref{4d.7}.

\noindent
2) Similar, the only difference is that the formula defining the
``order" is not first order and the length of the relevant sequences
may be infinite but still $\le |T|$ (see \cite[Ch.XIII]{Sh:c}).

\noindent
3) By parts (1),(2) we should consider only stable - not strongly
stable $T$.
Choose a class $\bold C$ of regular cardinals such that $\lambda
\in \bold C \Rightarrow (2^{< \lambda})^+ < \text{ Min}(\bold C
\backslash \lambda^+)$ and Min$(\bold C) > \mu$.  We iterate with full
support $\langle \bbP_\mu,\name{\bbQ}_\mu:\mu
\in \bold C\rangle$ with $\name{\bbQ}_\mu$ as in \ref{3c.4}. 
\end{PROOF}

\begin{claim}
\label{4d.17}
Assume $T \subseteq T_1,\lambda =
\lambda^\kappa$ is not necessary regular and $\kappa = \text{\rm
cf}(\kappa) < \kappa(T)$, e.g. $T$ is unstable.  \Then \, there are
{\rm EF}$^+_{\lambda \times \kappa,\lambda,\lambda^+}$-equivalent 
non-isomorphic models from {\rm PC}$(T_1,T)$ of cardinality $\lambda^+$.
\end{claim}

\begin{PROOF}{\ref{4d.17}}
As in \cite{HySh:529}, seeing the proof of \ref{4d.7}.
\end{PROOF}

\noindent
As said in the introduction by the old results (note: \ref{4d.23} is
on elementary classes and \ref{4d.27} on small enough pseudo
elementary classes).
\begin{conclusion}
\label{4d.23}
(ZFC) For first order countable
complete first order theory $T$ the following conditions are equivalent:
\mn
\begin{enumerate}
\item[$(A)$]  $T$ is superstable with NDOP and NOTOP
\sn
\item[$(B)_1$]   if $\lambda = \text{\rm cf}(\lambda) > |T|$ and $M_1,M_2
\in \text{ Mod}_T(\lambda)$ are 
$\bbL_{\infty,\lambda}(\tau_T)$-equivalent then $M_1,M_2$ are isomorphic
\sn
\item[$(B)_2$]  like (B)$_1$ for some $\lambda = 
\text{\rm cf}(\lambda) > |T|$
\sn
\item[$(C)$]  if $\lambda = \text{\rm cf}(\lambda) > |T|$ and $M_1,M_2
\in \text{ Mod}_T(\lambda^+)$ are EF$_{\omega,\lambda}$-equivalent
\then \, $M_1,M_2$ are isomorphic
\sn
\item[$(D)$]   for some regular $\lambda > |T|$, if
$M_1,M_2 \in \text{ Mod}_T(\lambda^+)$ are
EF$_{\lambda,\lambda^+}$-equivalent \then \, they are isomorphic.
\end{enumerate}
\end{conclusion}

\begin{PROOF}{\ref{4d.23}}
Clause (A), clause (B)$_1$, clause (B)$_2$ are
equivalent because: as proved in \cite[Ch.XIII,Th.1.11]{Sh:c}, we have $(A)
\Rightarrow (B)_1 \and (B)_2$ and the inverse
implication holds by \cite{Sh:220}.  Now by the definitions
trivially (B)$_1 \Rightarrow$ (C) $\Rightarrow$ (D).

Lastly, by \cite{HySh:529}, i.e. by \ref{4d.2} we 
have $\neg(A) \Rightarrow \neg(D)$,
i.e. (D) $\Rightarrow$ (A) so we have the circle.
\end{PROOF}

So \ref{4d.23} tells us what we know about Qustion $(A)_0$ of
\ref{1a.4}.  Similarly concerning $(B)_0$ of Question \ref{1a.4}.
\begin{conclusion}
\label{4d.27}
(ZFC)  For first order countable
complete first order theory $T$ and $\kappa \ge 2^{\aleph_0}$ 
the following conditions are equivalent:
\mn
\begin{enumerate}
\item[$(A)$]  $T$ is unsuperstable
\sn
\item[$(B)_\kappa$]  for every $\lambda > \kappa \ge |T|$ and
$(\kappa,T)$-candidate $\psi$ (see Definition \ref{0z.25}), and
ordinal $\alpha < \lambda$ satisfying $|\alpha|^+ = \lambda
\Rightarrow |\alpha \le |\alpha| \times \omega$, \underline{there are}
EF$_{\alpha,\lambda}$-equivalent non-isomorphic models $M_1,M_2 \in
\text{ PC}_{\tau(T)}(\psi)$ of cardinality $\lambda$
\sn
\item[$(C)$]   for some $\lambda > \kappa \ge |T|$, for no
$(\kappa,T)$-candidate $\psi$ is the class PC$_{\tau(T)}(\psi)$ 
categorical in $\lambda$.
\end{enumerate}
\end{conclusion}

\begin{PROOF}{\ref{4d.27}}
First, assume $T$ is superstable, so clause (A)
holds.  By the proofs of \cite[Ch.VI,\S4]{Sh:c} there is a
$(\kappa,T)$-candidate $\psi$, PC$_{\tau(T)}(\psi)$ is the class of
saturated models of $T$, (in details, if $n < \omega,\bar a \in
{}^n{\gC}$, tp$(b,\bar a,{\gC})$ is stationary, $q = \text{
tp}(\bar b,\emptyset,{\gC}),p=p(x,\bar y) = \text{ tp}(\langle a
\rangle \char 94 \bar b,\emptyset,{\gC})$ then let $\psi_{p,q}$ be such that $M
\models \psi_{p,q}$ \Iff \, for every $\bar b' \in {}^n M$ realizing
the type $q(\bar y)$, the function $c \mapsto F^M_{p,q}(c,\bar
b')$ is one-to-one and if $k < \omega,c_0,\dotsc,c_k \in M$ are
pairwise distinct then tp$_{\bbL(\tau(T))}
(F^M_{p,q}(c_k),\{F^M_{p,q}(c_0),\dotsc,F^M_{p,q}(c_{k-1})\} \cup 
\bar b',M)$ extends $p(x,\bar b')$ and does not fork over $M$.

Lastly, $\psi = \wedge\{\psi_{p,q}:p,q$ as above$\}$ so $\in 
\bbL_{\kappa^+,\omega}$.  So in the present case also
$(B)_\kappa,(C)_\kappa,(D)_\kappa$ holds. 

Second, assume $T$ is not superstable, so clause (A) holds and we
shall prove the rest.  Let $\psi$ be a $(\kappa,T)$-candidate.

By \ref{0z.29} and let there is $\Phi \in 
\Upsilon^{\omega_1{\text{-tr}}}_\kappa$ witnessing this hence
witnessing unsuperstability and now we can
use Theorem \ref{4d.2} quoted above.
\end{PROOF}
\newpage

\bibliographystyle{alphacolon}
\bibliography{lista,listb,listx,listf,liste,listz}

\end{document}